%% file: Matern_paper.tex
\documentclass[onefignum,onetabnum]{siamonline171218}

\usepackage[a4paper,margin=2.5cm,footskip=1cm]{geometry}
\usepackage{latexsym}
\usepackage{mathtools}
\usepackage{mathrsfs}
\usepackage{enumerate}
\usepackage{textcomp}
\usepackage[font=small]{caption}
\usepackage{tikz}
\usepackage{epstopdf}
\usepackage{tikz}
\usetikzlibrary{shapes,arrows,calc}

\usepackage{placeins}
\usepackage{graphicx}

\input{ex_shared}

\input{macros}

\newcommand{\eu}[1]{{\color{black}{#1}}}

\ifpdf
\hypersetup{
  pdftitle={Analysis of boundary effects on PDE-based sampling of Whittle-Mat\'ern random fields},
  pdfauthor={U. Khristenko, L. Scarabosio, P. Swierczynski, E. Ullmann, B. Wohlmuth}
}
\fi

\begin{document}

\maketitle

\begin{abstract}
We consider the generation of samples of a mean-zero Gaussian random field with Mat\'ern covariance function.
Every sample requires the solution of a differential equation with  Gaussian white noise forcing, formulated on a bounded computational domain.
This introduces unwanted boundary effects since the stochastic partial differential equation is originally posed on the whole $\bbR^d$, without boundary conditions. 
We use a window technique, whereby one embeds the computational domain into a larger domain, and postulates convenient boundary conditions on the extended domain.
To mitigate the pollution from the artificial boundary it has been suggested in numerical studies to choose a window size that is at least as large as the correlation length of the Mat\'ern field. 
We provide a rigorous analysis for the error in the covariance introduced by the window technique, for homogeneous Dirichlet, homogeneous Neumann, and periodic boundary conditions. 
We show that the error decays exponentially in the window size, independently of the type of boundary condition.
We conduct numerical experiments in 1D and 2D space, confirming our theoretical results.
\end{abstract}



\section{Introduction}

Gaussian random fields (GRFs) are an important building block in computational statistics and uncertainty quantification, where they are typically used as inputs to complex models of physical and technical processes.
Common applications arise for instance in 
astrophysics \cite{Bardeen1986},
biology \cite{TW2007},
geophysics \cite{BT2013, Isaac2015},
geotechnical engineering \cite{GHF2009}, 
hydrology \cite{SGC2006},
image processing \cite{Cohen1991}, and
meteorology \cite{BL2011, LRL}.
The popularity of GRFs can be attributed to the fact that they are simple yet flexible models of spatial variability.
A GRF is completely characterized by the mean function and covariance operator.
As a second order random field, a GRF admits an expansion in eigenpairs of the covariance operator, this is known as Karhunen-Lo\`eve (KL) expansion \cite{GS1991, LPS2014}.
Moreover, depending on the regularity of the covariance operator, it is possible to produce realizations with different smoothness properties. 
However, the efficient generation of these realizations is a serious computational bottleneck for several reasons. 
A random function can be considered as a random variable taking values in an infinite-dimensional function space.
Alternatively, it can be considered as a collection of random variables defined on a suitable probability space, indexed by their spatial position.
In either case, a joint discretization of a function space coupled with a probability space is necessary in practice.
This may lead to \textit{high-dimensional} approximation spaces with prohibitive computational costs and memory requirements.
In addition, many commonly used covariance operators are \textit{nonlocal}, and couple the random variations within the spatial domain of interest.
This means that in a realization the function value at a specific location depends on contributions from all other locations.

A simple strategy to sample a GRF at a fixed number of spatial grid points uses a factorization of the covariance matrix associated with these points.
However, due to the nonlocality of many covariance functions, including the Mat\'ern covariance which we discuss below, the covariance matrix is a large, dense matrix, which is expensive to factorize.
The computational cost can be reduced by efficient approximations of the covariance matrix based on low-rank approaches (see e.g. \cite{BennerQiuStoll2018, DElia2013, Harbrechtetal2012, Harbrechtetal2015}) or hierarchical matrices (see e.g. \cite{ChenStein2017, Feischl2018}).
Alternatively, we can employ circulant embedding \cite{ChanWood1997, DN, Graham18}, where the factorization is performed by the Fast Fourier Transform (FFT). 
This is very efficient but requires simple geometries together with uniform, structured grids, and stationary covariance operators, and might limit scalability in a massively  parallel environment \cite{Dz}. 
In addition, the circulant matrix might have negative eigenvalues which have to be discarded and thus introduce an error in the covariance representation \cite[\S 6.5]{LPS2014}.
Alternatively, the embedding can be chosen large enough to guarantee non-negative eigenvalues; this introduces additional costs \cite{Graham18}.

Another option to generate realizations of a GRF is through its KL expansion \cite{Betzetal2014, Eiermann2007}. 
Unfortunately, it usually requires the solution of an integral eigenproblem with a large, dense matrix caused again by a nonlocal covariance operator.
This can be made more efficient by similar techniques used for the covariance matrix factorization (see e.g. \cite{Khoromskij2009, Kressner2015, SchwabTodor2006}).
Alternative ideas for an efficient computation of the KL expansion are based on domain decomposition \cite{Rizzi2017}, randomized linear algebra \cite{SLK}, or domain modification \cite{PraneshGhosh2015}.
However, if the eigenvalues of the covariance operator decay slowly, which is the case for small correlation lengths or fields with very low spatial regularity, then many terms in the KL expansion are needed to achieve a good approximation for the GRF.
This leads to high computational costs for each sample, and large memory requirements for storing the KL eigenvalues and eigenfunctions. 


A third alternative is to compute each sample of a GRF as a solution to a fractional elliptic stochastic partial differential equation (SPDE) with white noise forcing.
For GRFs with Mat\'ern covariance function the connection to SPDEs has been noticed by Whittle in \cite{Whittle54, Whittle63}.
This approach allows for mesh flexibility and is more efficient than the KL expansion for slowly decaying eigenvalues.
Moreover, we can build on well-established discretizations for partial differential equations (PDEs) and associated fast, scalable multigrid solvers \cite{Osb17, Osb18}.
However, the efficient implementation and rigorous error analysis of PDE-based sampling is not straightforward since the SPDE is posed on the full space $\mathbb{R}^d$ \textit{without} boundary conditions.
For practical simulations a truncation to some bounded domain of interest $D \subset \mathbb{R}^d$ is necessary.
This also requires boundary conditions for the samples which are not known in general.
Finally, the \eu{SPDE} requires a suitable discretization.
We discuss these points in more detail in the following.

Lindgren et al. \cite{LRL} introduced a finite element discretization of the fractional elliptic SPDE on bounded domains augmented with artificial Neumann boundary conditions. 
This approach has become standard in the recent literature \cite{Croci18,Dz,Osb17,Osb18,Roininen14} and is termed \textit{window technique}.
Here we embed the domain $D$ into a larger domain and solve the SPDE on the latter, using homogeneous Neumann or Dirichlet boundary conditions. 
In this way, the boundary effects coming from the domain truncation are negligible in $D$ if the window is large enough. 
In \cite{LH15, Roininen14} an empirical rule for the window size is suggested, where it is observed that the window boundary should be at least as far away from $D$ as the correlation length of the Mat\'ern field.
However, a precise error analysis has not been carried out to date.
In this paper, we provide an analysis of the boundary effects on the covariance structure of the samples of the Mat\'ern field, thereby closing a gap in the literature.
We study the error in the covariance function depending on the window size and for different types of boundary conditions, namely, homogeneous Dirichlet, homogeneous Neumann and periodic boundary conditions.
We show that the domain truncation introduces an \textit{aliasing effect}, and that the covariance error decays exponentially in the window size, independently of the type of boundary condition considered.
Moreover, we provide numerical results for a specific choice of Robin boundary conditions.

Our choice of the boundary conditions is guided by the literature, and by practical considerations.
Homogeneous Dirichlet \cite{Croci18, Roininen14} resp. homogeneous Neumann conditions \cite{Dz, Osb17, Osb18, Roininen14} are easy to implement in finite element frameworks and for FFT on tensor product domains, and they do not require additional parameter tuning. 
Although periodic boundary conditions are not popular in this context, we include these here since they mimic the stationarity of the exact Mat\'ern samples on the full space. 
Moreover, they give rise to analytic expressions for the covariance error. 
Robin boundary conditions for the SPDE on the bounded domain have been considered in \cite{Daon2018, Roininen14}. 
Robin conditions involve a coefficient which has to be tuned. 
Furthermore, the \eu{structure} of eigenvalues and eigenfunctions for this case requires a different analysis for the error compared to the previously mentioned boundary conditions.
\eu{Finally, Robin boundary conditions cannot be trivially applied if the SPDE is solved with the FFT.
For these reasons, we only} consider a specific choice of coefficient for Robin boundary conditions, motivated by the results in \cite{Daon2018} and some physical considerations, and show numerically its effects on the covariance of the sampled field.

We mention that the window technique is not the only option to address the unwanted boundary effects in PDE-based sampling.
Alternatively, one can try and find an accurate approximation to the \textit{exact} boundary condition on the domain $D$ without a window.
Unfortunately, the exact boundary conditions are only known for very special parameter configurations in the Mat\'ern covariance.
Daon and Stadler \cite{Daon2018} show that a specific Robin condition is exact in 1D space for the exponential covariance, however, the general case is an open research question. 
They suggest to optimize a varying Robin coefficient, or to rescale the covariance operator, which again only mitigates the boundary effect.

\eu{
Finally, we observe that PDE-based sampling requires a discretization of the fractional elliptic SPDE; this is often performed together with the white noise discretization. 
Recent studies on the finite element discretization and error analysis of semilinear and fractional elliptic SPDEs with white noise forcing can be found in \cite{Zhang16} and \cite{BK17, BKK18, BKK17}, respectively.
For the spatial discretization, classical piecewise linear \cite{BKK18, BKK17, Croci18, LRL} or mixed finite elements \cite{Osb17, Osb18} have been employed.
However, }if the white noise is discretized by finite elements, then the load vector follows a Gaussian distribution with the finite element mass matrix as covariance matrix. 
Sampling from it entails therefore the factorization of the mass matrix.
In finite element frameworks, the mass matrix is usually sparse, however, its factorization cost will not scale optimally in the number of degrees of freedom.
This issue can be resolved by mass lumping \cite{LRL}, or by approximating the white noise by a suitable piecewise constant random function \cite{ANZ, Dz, DuZhang03}.
Recently, Croci et al. \cite{Croci18} proposed a white-noise sampling \eu{with optimal, linear complexity}, where the finite element mass matrix is factorized exactly using small element mass matrices.
We envision that this will further enhance the popularity of PDE-based sampling of Mat\'ern fields. 
Our analysis provides a missing piece of information on another ingredient, the choice of boundary conditions, for this attractive and efficient sampling framework.

The paper is organized as follows. 
In \Cref{sect:formulation}, we review PDE-based sampling, and introduce the problem formulation on the full and bounded spatial domain, respectively. 
In \Cref{sect:folded}, we study the aliasing effect on the covariance introduced by the domain truncation. 
The results of \Cref{sect:folded} are used in \Cref{sect:mainthm} to derive the main contribution of this paper, that is an error bound for the Mat\'ern covariance in terms of the window size. 
\eu{In \Cref{sect:ext} we discuss the extension of the error analysis for \textit{anisotropic} Mat\'ern covariances, and provide a motivation for the use of Robin boundary conditions.}
In \Cref{sect:numexp}, we show numerical experiments which confirm the theory of the previous sections. 
In the same section, we provide numerical evidence that a specific choice of Robin boundary conditions, possibly not optimal but motivated by previous work \cite{Daon2018} and analogies with the Helmholtz equation, has better approximation properties in the covariance than the other boundary conditions considered.
Finally, \Cref{sect:concl} presents some concluding remarks.

\section{PDE-based sampling}\label{sect:formulation}

Let $D\subset\mathbb{R}^d$ be an open, simply connected, bounded domain, for $d=1,2,3$, and let $(\Omega,\mathcal{A},\mathbb{P})$ be a probability space (with $\Omega$ the set of events, $\mathcal{A}$ a $\sigma$-algebra and $\mathbb{P}$ a probability measure). 
We consider the task of sampling a Gaussian field $u=u(\omega,\bfx)$, $\omega\in\Omega$ and $\bfx\in D$, with zero mean and Mat\'ern covariance function \cite{Guttorp2006, Stein1999}. 
That is, for every $\bfx,\bfy\in\mathbb{R}^d$, the covariance function of $u$ is given by
\begin{equation}\label{eq:matcov}
\mathcal{C}(\bfx,\bfy) 
= 
\sigma^2 \mathcal{M}_\nu(\kappa\Vert \bfx-\bfy\Vert_2),\quad \kappa = \frac{\sqrt{2\nu}}{\rho},
\end{equation}
where $\sigma^2>0$ is the marginal variance, $\kappa$ is a scaling parameter, $\rho$ is the 
correlation length, and the unit Mat\'ern function $\mathcal{M}_\nu(x)$ is defined as
\begin{equation}\label{eq:matcov2}
\mathcal{M}_\nu(x) 
= 
\frac{x^\nu \mathcal{K}_\nu(x)}{2^{\nu-1}\Gamma(\nu)},
\end{equation}
for $x>0$, with $\mathcal{K}_{\nu}$ the modified Bessel function of the second  kind and order $\nu>0$.
The parameter $\nu$ determines the mean square differentiability of $u$ \cite[\S 4.2]{RW}. 
The Mat\'ern covariance kernel is then fully characterized by the three deterministic parameters $\sigma^2$, $\rho$ and $\nu$, which vary independently. 
It is clear from \eqref{eq:matcov} that the Mat\'ern covariance is invariant under translations and rotations, hence the Gaussian field $u$ is {isotropic} (see e.g. \cite[p. 33]{Adler1981}).
This is in fact a special case of a {stationary} random field \cite[p. 24]{Adler1981}.


\subsection{Formulation on the full spatial domain}
In \cite{Whittle54, Whittle63} Whittle showed that a zero mean Gaussian field with covariance function \eqref{eq:matcov} is the solution to the SPDE
\begin{equation}\label{eq:spde}
\left(\mathcal{I}-\kappa^{-2}\Delta\right)^{\frac{\alpha}{2}} u(\omega,\bfx) = \eta\dot{W}(\omega,\bfx),\quad \bfx\in\mathbb{R}^d,\quad \text{for }\mathbb{P}\text{-a.e. } \omega\in\Omega,
\end{equation}
where $\mathcal{I}$ is the identity operator, $\dot{W}$ denotes Gaussian white noise as defined in \cite[Def. 6 p.~448]{LRL}, $\alpha = \nu + \frac{d}{2}$, and $\eta$ is a normalization constant set to
\begin{equation}
\eta^2 = \sigma^2\frac{\left(4\pi\right)^{d/2}\Gamma(\nu+d/2)}{\kappa^d\Gamma(\nu)},
\end{equation}
guaranteeing that the marginal variance of $u$ at each point in the domain is equal to $\sigma^2$. 
Equation \eqref{eq:spde} has to be understood in the sense of distributions. 
The SPDE is posed on $\mathbb{R}^d$ \textit{without} boundary conditions.
However, for simulations, a truncation of $\mathbb{R}^d$ to some bounded domain of interest $D$ is necessary;
this also requires boundary conditions for $u$.
Here we employ the \textit{window technique} \cite{Croci18,Dz,LRL,Osb17,Osb18,Roininen14}, where we embed the domain $D$ into a larger domain and solve \eqref{eq:spde} on the latter, using inexact, artificial boundary conditions for $u$. 
We will see that by truncating the domain an aliasing effect occurs, referred to in \cite{LRL} as \textit{folded covariance}.
However, the intuition is that the boundary effects coming from the domain truncation are negligible in $D$ if the window is large enough.
In \cite{LH15,Roininen14} the empirical rule of taking a window with boundaries at distance at least $\rho$ from $D$ is suggested. 
To investigate this, we consider the SPDE \eqref{eq:spde} on a bounded domain, and derive explicit expressions for the folded covariance.


\subsection{Formulation on a bounded domain}
\begin{figure}
\centering
 \begin{tikzpicture}
            \filldraw[white!70!cyan,draw=blue]  plot[smooth, tension=0.8] coordinates {(0.1,3) (0.5,3.5) (1.3,3.5) (2,3.3) (2.5,3.5) (3,2.8)  (3.98,0.5) (1.5,0.06) (0,0.8) (0.1,3)};
            \node at (1.5,1.5) {$D$};
            \node at (5.3,3.5) {$D_{ext}$};
            \draw (-0.95,-0.8) rectangle (4.85,4.7);
            \draw[->,thick] (-0.95,-0.8)--(-0.95,5.45);
            \draw[->,thick] (-0.95,-0.8)--(5.55,-0.8);
            \coordinate (L) at ($(4.85,-0.8)$) {};
      		\draw ($(L)+(0,3pt)$) -- ($(L)-(0,3pt)$);
      		\node at ($(L)-(0,2.5ex)$) {$L$};
      		\coordinate (dr) at ($(4.05,-0.8)$) {};
      		\draw ($(dr)+(0,3pt)$) -- ($(dr)-(0,3pt)$);
      		\node at ($(dr)-(0,2.5ex)$) {$L-\tfrac{\delta}{2}$};
      		\coordinate (dl) at ($(-0.15,-0.8)$) {};
      		\draw ($(dl)+(0,3pt)$) -- ($(dl)-(0,3pt)$);
      		\node at ($(dl)-(0,2.5ex)$) {$\tfrac{\delta}{2}$};
      		\node at (-0.95,-1) {$0$};
      		\coordinate (Ll) at ($(-0.95,4.7)$) {};
      		\draw ($(Ll)-(0.15,0pt)$) -- ($(Ll)+(0.15,0pt)$);
      		\node at ($(Ll)-(0.5,0ex)$) {$L$};
      		\coordinate (du) at ($(-0.95,4)$) {};
      		\draw ($(du)-(0.15,0pt)$) -- ($(du)+(0.15,0pt)$);
      		\node at ($(du)-(0.82,0ex)$) {$L-\tfrac{\delta}{2}$};
      		\coordinate (dd) at ($(-0.95,0)$) {};
      		\draw ($(dd)-(0.15,0pt)$) -- ($(dd)+(0.15,0pt)$);
      		\node at ($(dd)-(0.4,0ex)$) {$\tfrac{\delta}{2}$};
      		\draw[dashed] (dr)--(4.05,4.7);
      		\draw[dashed] (dl)--(-0.15,4.7);
      		\draw[dashed] (du)--(4.85,4);
      		\draw[dashed] (dd)--(4.85,0);
\end{tikzpicture}\caption{Embedding of the domain $D$ in the extended domain $D_{ext}=(0,L)^2$.}\label{fig:window}
\end{figure}
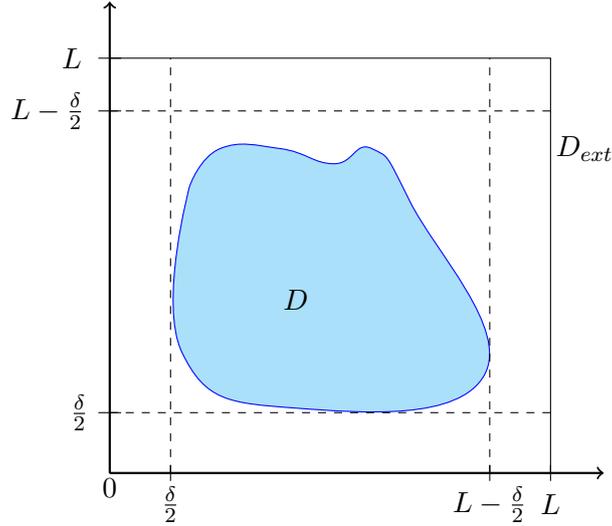

We consider a coordinate system such that $\inf_{\bfx\in D} x_i = \frac{\delta}{2}$, for $i=1,\ldots,d$ and $\delta>0$, and where $x_i$, $i=1,\ldots,d$, denotes the $i$-th coordinate of $\bfx$. 
We define the extended domain in which the SPDE in \eqref{eq:spde} has to be solved as $D_{ext}:=(0,L)^d$, for $L=\delta+\ell$, and $\ell=\sup_{\bfx,\bfy\in D}\lVert \bfx-\bfy\lVert_{\infty}$. 
The setting is depicted in \cref{fig:window}. 

We consider then the solution of \eqref{eq:spde} on $D_{ext}$ and denote it by $u_L$. 
The reason for choosing a tensor-product domain for $D_{ext}$ is motivated \eu{by the fact} that this allows for easy application of the FFT to obtain samples for $u_L$.
For the differential operator, we use the spectral definition of the fractional Laplacian, in agreement with its use in the statistics literature \cite{LRL}. 
We address the approximation of the Mat\'ern covariance when imposing either one of the following boundary conditions: \smallskip
\begin{enumerate}
\item[($D$)] homogeneous Dirichlet:\quad $\left.(-\Delta)^j u_L\right|_{\partial D_{ext}}=0$, $j=0,\ldots,\lfloor \frac{\alpha-1}{2}\rfloor$,\smallskip
\item[($N$)] homogeneous Neumann:\quad $\left.\dfrac{\partial}{\partial \bfn}(-\Delta)^ju_L\right|_{\partial D_{ext}}=0$, $j=0,\ldots,\lfloor \frac{\alpha-1}{2}\rfloor$,\smallskip
\item[($P$)] periodic:
\quad $\left.(-\Delta)^j u_L\right|_{x_i=0}=\left.(-\Delta)^j u_L\right|_{x_i=L}$ and $\left.\dfrac{\partial}{\partial \bfn}(-\Delta)^j u_L\right|_{x_i=0}=\left.\dfrac{\partial}{\partial \bfn}(-\Delta)^j u_L\right|_{x_i=L}$, $j=0,\ldots,\lfloor \frac{\alpha-1}{2}\rfloor$ and $i=1,\ldots,d$,\smallskip
\end{enumerate}
for $\mathbb{P}$-a.e. $\omega\in\Omega$, where $\bfn$ denotes the outward unit normal to $\partial D_{ext}$.

Note that by construction the average of the random field $u_L$ on the truncated domain $D_{ext}$ is zero. 
Indeed, for any boundary condition the solution can be written as the It\^{o} integral \cite[Appendix B]{Holden96}
\begin{equation}\label{eq:ulgreen}
u_L(\bfx)=\int_{D_{ext}} G_{\frac{\alpha}{2}}(\bfx,\bfy)\;\dd \dot{W}(\bfy), \quad \bfx\in D_{ext},
\end{equation}
where $G_{\frac{\alpha}{2}}$ is the corresponding Green's function. Equation \eqref{eq:ulgreen} is well-defined and holds almost surely due to \cite[Lemma 6.25]{Stuart}.
Hence $\bbE[u_L]=0$ by definition of the integral. 

Let $\mathcal{C}^{L}_{\ast}$, for $\ast\in\left\{D,N,P\right\}$, denote the covariance function of $u_L$ when using Dirichlet, Neumann or periodic boundary conditions, respectively. 
This is the Green's function of $\left(\mathcal{I}-\kappa^{-2}\Delta\right)^{\alpha}$ with corresponding boundary conditions \cite{Daon2018}.
To study the convergence for the second moment of the random field $u$, we study the error
\begin{equation}\label{eq:linferr}
\lVert \mathcal{C}-\mathcal{C}^{L}_{\ast}\rVert_{\infty}:=\sup_{\bfx,\bfy\in D}|\mathcal{C}(\bfx,\bfy)-\mathcal{C}^{L}_{\ast}(\bfx,\bfy)|,
\end{equation}
for $\ast\in\left\{D,N,P\right\}$.

\subsection{Covariance of the solution on the bounded domain}\label{sect:folded}

In this section, we extend the result in \cite[Thm. 1]{LRL} to Dirichlet and periodic boundary conditions, and to the spatial dimension $d>1$. 
A central fact that is used in the proofs in this and in the next section is the following.
The covariance functions $\mathcal{C}^{L}_{\ast}$, for $\ast\in\left\{D,N,P\right\}$, admit the expansion (see e.g. \cite[Eq. 35]{LRL})
\begin{equation}\label{eq:covexp}
\mathcal{C}^{L}_{\ast}(\bfx,\bfy) = \eta^2\sum_{\bfk\in\mathcal{N}^d} \lambda_{\bfk}^{-\alpha}w_{\bfk}(\bfx)w_{\bfk}(\bfy),
\end{equation}
where $\mathcal{N}=\bbN$ or $\mathcal{N}=\bbN_0:=\bbN\cup\left\{0\right\}$ according to the boundary conditions, $\left\{w_{\bfk}\right\}_{\bfk\in\mathcal{N}^d}$ are the eigenfunctions of the operator $(\mathcal{I}-\kappa^{-2}\Delta)$ on $D_{ext}$ and $\left\{\lambda_{\bfk}\right\}_{k\in\mathcal{N}^d}$ are the associated eigenvalues.  
In particular, the eigenpairs $(\lambda_{\bfk},w_{\bfk})$ in \eqref{eq:covexp} are  
\begin{equation}\label{eq:eigdir}
w_{\bfk}= \prod_{i=1}^d a_{k_i}\sin\left(\frac{\pi k_i x_i}{L}\right),\quad 
\lambda_{\bfk}=1+\left(\frac{\pi}{\kappa L}\right)^2\lVert\bfk\rVert_2^2
\end{equation}
and $\mathcal{N}=\bbN$ for Dirichlet boundary conditions,
\begin{equation}\label{eq:eigneu}
w_{\bfk}= \prod_{i=1}^d a_{k_i}\cos\left(\frac{\pi k_i x_i}{L}\right),\quad 
\lambda_{\bfk}=1+\left(\frac{\pi}{\kappa L}\right)^2\lVert\bfk\rVert_2^2
\end{equation}
and $\mathcal{N}=\bbN_0$ for Neumann boundary conditions, and, denoting $\img=\sqrt{-1}$,
\begin{equation}\label{eq:eigper}
w_{\bfk}= \frac{1}{L^{\frac{d}{2}}}e^{\frac{2\pi}{L}\img\bfk\cdot\bfx},\quad 
\lambda_{\bfk}=1+\left(\frac{2\pi}{\kappa L}\right)^2 \lVert\bfk\rVert_2^2
\end{equation}
and $\mathcal{N}=\bbN_0$ for periodic boundary conditions. 
The normalizing constants are $a_{k_i}=\sqrt{\frac{1}{L}}$ for $k_i=0$ and $a_{k_i}=\sqrt{\frac{2}{L}}$ otherwise. 
Note that, by the spectral definition of the fractional Laplacian, the operators $(\mathcal{I}-\kappa^{-2}\Delta)$ and $(\mathcal{I}-\kappa^{-2}\Delta)^{\frac{\alpha}{2}}$, $\alpha>\frac{d}{2}$, have the exact same eigenfunctions \cite{lectnotes}. 
Analogously to \cite[Thm. 1]{LRL}, we obtain the following result.

\begin{theorem}\label{thm:folded}
The covariance function of the solution $u_L$ to \eqref{eq:spde} when using, respectively, periodic, homogeneous Neumann and homogeneous Dirichlet boundary conditions on $\partial D_{ext}$ is given by
\begin{align}
		\mathcal{C}_{P}^L(\bfx,\bfy) &= \sum_{\bfk\in\mathbb{Z}^d}\mathcal{C}(\bfx+L\bfk,\bfy), \label{eq:cov_P} \\
		\mathcal{C}_{N}^{L}(\bfx,\bfy) & = \sum_{\boldsymbol{\varepsilon}\in B^d} \mathcal{C}_{P}^{2L}(\bfx, \boldsymbol{\varepsilon}.\bfy)
		= \sum_{\bfeps\in B^d}\sum_{\bfk\in\bbZ^d} \Matern\left(\bfx + 2L\bfk,\bfeps.\bfy\right) ,  \label{eq:cov_N} \\
		\mathcal{C}_{D}^{L}(\bfx,\bfy) & =  \sum_{\boldsymbol{\varepsilon}\in B^d}\left(\prod_{i=1}^{d}\varepsilon_i\right)\mathcal{C}_{P}^{2L}(\bfx, \boldsymbol{\varepsilon}.\bfy) = \sum_{\boldsymbol{\varepsilon}\in B^d}\left(\prod_{i=1}^{d}\varepsilon_i\right)\sum_{\bfk\in\bbZ^d} \Matern\left(\bfx + 2L\bfk,\bfeps.\bfy\right),  \label{eq:cov_D}
		\end{align}
		where $B=\left\{-1,1\right\}$ and $\bfx.\bfy=(x_1 y_1,\ldots,x_d y_d)^\top$ denotes the element-wise product.
\end{theorem}
\begin{proof}
For fixed $\bfx, \bfy\in D$, the folding effect introduced by the domain truncation is illustrated in \cref{fig:folding}, for all three boundary conditions.
In the case of periodic boundary conditions, according to \eqref{eq:covexp} and \eqref{eq:eigper}, we have
	\begin{equation}
	\cov_{P}^L(\bfx,\bfy)=\frac{\eta^2}{L^{d}}\sum_{\bfk\in\mathbb{Z}^d}\left(1 + \left(\frac{2\pi}{\kappa L}\right)^2 \lVert\bfk\rVert_2^2 \right)^{-\alpha}\, e^{\frac{2\pi \img}{L} \bfk \cdot(\bfx-\bfy)}. \label{eq:FourierCovP}
	\end{equation}
Given the Fourier representation of the Matern kernel~\cite{Roininen14}
\begin{equation*}
	\cov(\bfx,\bfy) = \frac{\eta^2}{(2\pi)^d}\int_{\bbR^d}  \left(1 + \kappa^{-2}\lVert\bfxi\rVert_2^2 \right)^{-\alpha}\, e^{\img \bfxi \cdot(\bfx-\bfy)}  \dd\bfxi,
\end{equation*}
a direct application of the Poisson summation formula~\cite{helson1983harmonic,stein2016introduction} to~(\ref{eq:FourierCovP}) yields
\begin{equation*}
\mathcal{C}_{P}^L(\bfx,\bfy) = \sum_{\bfk\in\bbZ^d} \mathcal{C}(\bfx + L\bfk, \bfy ).
\end{equation*}
For Dirichlet boundary conditions it holds
\begin{align}
\cov_D^L(\bfx,\bfy) & = \frac{\eta^2}{L^d}\sum_{\bfk\in\bbZ^d}\lambda_{\bfk}^{-\alpha}\prod_{i=1}^d \sin\left(\frac{\pi k_i x_i}{L}\right)\sin\left(\frac{\pi k_i y_i}{L}\right)\nonumber\\
&= \frac{\eta^2}{L^d}\sum_{\bfk\in\bbZ^d}\lambda_{\bfk}^{-\alpha}\cdot  \frac{1}{4^d}\sum_{\bfeps_1,\bfeps_2\in B^d}\left(\prod_{j=1}^{d}\bfeps_{1,j}\bfeps_{2,j}\right) e^{\img\frac{\pi}{L}\bfk\cdot\left(\bfeps_1.\bfx-\bfeps_2.\bfy\right)}\nonumber\\
&= \sum_{\bfeps\in B^d}\left(\prod_{j=1}^{d}\bfeps_{j}\right)\frac{\eta^2}{(2L)^d}\sum_{\bfk\in\bbZ^d} \lambda_{\bfk}^{-\alpha}\, e^{\img\frac{2\pi}{2L}\bfk\cdot\left(\bfx-\bfeps.\bfy\right)}\nonumber\\
&= \sum_{\bfeps\in B^d}\left(\prod_{j=1}^{d}\bfeps_{j}\right)\cov_{P}^{2L}(\bfx,\bfeps.\bfy), \label{eq:dirsum}
\end{align}
where we have used Euler's formulas in the second line and the eigenvalues are given by \eqref{eq:eigdir}.
Similar computations lead to \eqref{eq:cov_N} for Neumann boundary conditions:
\begin{align*}
\cov_N^L(\bfx,\bfy) & = \frac{\eta^2}{L^d}\sum_{\bfk\in\bbZ^d}\lambda_{\bfk}^{-\alpha}\prod_{i=1}^d \cos\left(\frac{\pi k_i x_i}{L}\right)\cos\left(\frac{\pi k_i y_i}{L}\right)\\
&= \frac{\eta^2}{L^d}\sum_{\bfk\in\bbZ^d}\lambda_{\bfk}^{-\alpha}\cdot  \frac{1}{4^d}\sum_{\bfeps_1,\bfeps_2\in B^d} e^{\img\frac{\pi}{L}\bfk\cdot\left(\bfeps_1.\bfx - \bfeps_2.\bfy\right)}\\
&= \sum_{\bfeps\in B^d}\cov_{P}^{2L}(\bfx,\bfeps.\bfy).\\
\end{align*}
with $\left\{\lambda_{\bfk}\right\}_{\bfk\in\bbZ^d}$ as from \eqref{eq:eigneu}.
\end{proof}

\begin{figure}
\centering
 \begin{tikzpicture}[scale=1.5]
 \draw[->,thick] (-2.4,0)--(2.8,0);
 \draw[->,thick] (0,-1.4)--(0,2.8);
 \node at (2.6,-0.2) {$x_1$};
 \node at (0.2,2.6) {$x_2$};
 \coordinate (x) at ($(0.3,0.3)$) {}; 
 \filldraw (x) circle (0.03cm);
 \node at ($(x)-(0.1,0.1)$) {$\bfx$};
 \coordinate (y) at ($(0.4,0.4)$) {};
 \filldraw (y) circle (0.03cm);
 \node at ($(y)+(0.1,0.1)$) {$\bfy$};
 \node at (0.1,1.1) {$L$};
 \node at (1.1,-0.1) {$L$};
 \draw[thin] (-2.4,1)--(2.4,1);
 \draw[thin] (-2.4,2)--(2.4,2);
 \draw[thin] (-2.4,-1)--(2.4,-1);
 \draw[thin] (-2,-1.4)--(-2,2.4);
 \draw[thin] (-1,-1.4)--(-1,2.4);
  \draw[thin] (2,-1.4)--(2,2.4);
 \draw[thin] (1,-1.4)--(1,2.4);
 \filldraw[fill=black!10!green,draw=black!10!green] ($(y)+(1,0)$) circle (0.03cm);
 \filldraw[fill=black!10!green,draw=black!10!green] ($(y)+(1,1)$) circle (0.03cm);
 \filldraw[fill=black!10!green,draw=black!10!green] ($(y)+(0,1)$) circle (0.03cm);
 \filldraw[fill=black!10!green,draw=black!10!green] ($(y)+(-1,0)$) circle (0.03cm);
 \filldraw[fill=black!10!green,draw=black!10!green] ($(y)+(-1,-1)$) circle (0.03cm);
 \filldraw[fill=black!10!green,draw=black!10!green] ($(y)+(0,-1)$) circle (0.03cm);
 \filldraw[fill=black!10!green,draw=black!10!green] ($(y)+(-2,0)$) circle (0.045cm);
 \filldraw[fill=black!10!green,draw=black!10!green] ($(y)+(-2,-1)$) circle (0.03cm);
 \filldraw[fill=black!10!green,draw=black!10!green] ($(y)+(-1,1)$) circle (0.03cm);
 \filldraw[fill=black!10!green,draw=black!10!green] ($(y)+(-2,1)$) circle (0.03cm);
 \filldraw[fill=black!10!green,draw=black!10!green] ($(y)+(1,-1)$) circle (0.03cm);
 \coordinate (y1) at ($(-0.4,0.4)$) {};
 \coordinate (y2) at ($(-0.4,-0.4)$) {};
 \coordinate (y3) at ($(0.4,-0.4)$) {};
 \filldraw[fill=red!60!blue,draw=red!60!blue] ($(y1)$) circle (0.03cm);
  \node at($(y1)+(0.1,0.1)$) {\small\color{blue}{$+$}};
 \node at($(y1)+(0.1,-0.1)$) {\small\color{red}{$-$}};
 \filldraw[fill=red!60!blue,draw=red!60!blue] ($(y2)$) circle (0.03cm);
  \node at($(y2)+(0.1,0.1)$) {\small\color{blue}{$+$}};
 \node at($(y2)+(0.1,-0.1)$) {\small\color{red}{$+$}};
 \filldraw[fill=red!60!blue,draw=red!60!blue] ($(y3)$) circle (0.03cm);
 \node at($(y3)+(0.1,0.1)$) {\small\color{blue}{$+$}};
 \node at($(y3)+(0.1,-0.1)$) {\small\color{red}{$-$}};
 \filldraw[fill=red!60!blue,draw=red!60!blue] ($(y1)+(2,0)$) circle (0.03cm);
 \node at($(y1)+(2.1,0.1)$) {\small\color{blue}{$+$}};
 \node at($(y1)+(2.1,-0.1)$) {\small\color{red}{$-$}};
 \filldraw[fill=red!60!blue,draw=red!60!blue] ($(y2)+(2,0)$) circle (0.03cm);
  \node at($(y2)+(2.1,0.1)$) {\small\color{blue}{$+$}};
 \node at($(y2)+(2.1,-0.1)$) {\small\color{red}{$+$}};
  \filldraw[fill=red!60!blue,draw=red!60!blue] ($(y2)+(0,2)$) circle (0.03cm);
   \node at($(y2)+(0.1,2.1)$) {\small\color{blue}{$+$}};
 \node at($(y2)+(0.1,1.9)$) {\small\color{red}{$+$}};
 \filldraw[fill=red!60!blue,draw=red!60!blue] ($(y3)+(0,2)$) circle (0.03cm);
  \node at($(y3)+(0.1,2.1)$) {\small\color{blue}{$+$}};
 \node at($(y3)+(0.1,1.9)$) {\small\color{red}{$-$}};
 \filldraw[fill=red!60!blue,draw=red!60!blue] ($(y2)+(2,2)$) circle (0.03cm);
  \node at($(y2)+(2.1,2.1)$) {\small\color{blue}{$+$}};
 \node at($(y2)+(2.1,1.9)$) {\small\color{red}{$+$}};
 \filldraw[fill=red!60!blue,draw=red!60!blue] ($(y)-(2,0)$) circle (0.03cm);
  \node at($(y)+(-1.9,0.1)$) {\small\color{blue}{$+$}};
 \node at($(y)+(-1.9,-0.1)$) {\small\color{red}{$+$}};
 \filldraw[fill=red!60!blue,draw=red!60!blue] ($(y3)-(2,0)$) circle (0.03cm);
  \node at($(y3)+(-1.9,0.1)$) {\small\color{blue}{$+$}};
 \node at($(y3)+(-1.9,-0.1)$) {\small\color{red}{$-$}};
 \filldraw[fill=red!60!blue,draw=red!60!blue] ($(y3)+(-2,2)$) circle (0.03cm);
  \node at($(y3)+(-1.9,2.1)$) {\small\color{blue}{$+$}};
 \node at($(y3)+(-1.9,1.9)$) {\small\color{red}{$-$}};
 \end{tikzpicture}\caption{Folding effect due to domain truncation, for $d=2$. For fixed $\bfx, \bfy\in D$, the covariance $\cov_P^L(\bfx,\bfy)$ when using periodic boundary conditions is given by $\cov(\bfx,\bfy)$ plus the Mat\'ern covariance between $\bfx$ and all the green dots ($\bfy$ translated by $\bfk L$, $\bfk\in\bbZ^d$). For Neumann and Dirichlet boundary conditions, the covariance is given by $\cov(\bfx,\bfy)$ and the contributions given by the violet dots ($\bfeps.\bfy$ translated by $2\bfk L$, $\bfk\in\bbZ^d$, $\bfeps\in \left\{-1,1\right\}^d$). In the Neumann case, all these contributions are positive (blue signs), while the signs alternate for the Dirichlet case (see red signs in the figure).}\label{fig:folding}
\end{figure}
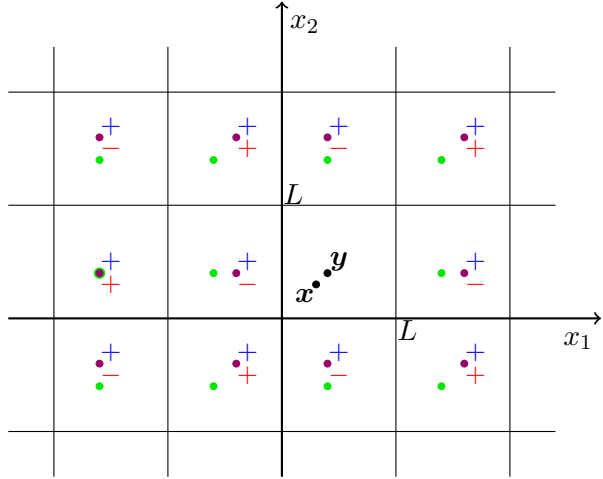
\begin{remark}
Recall that the Mat\'ern kernel on the full domain $\mathbb{R}^d$ is invariant under translations and rotations.
However, the expressions in \eqref{eq:cov_P}--\eqref{eq:cov_D} clearly tell us that this is not the case on the truncated domain $D_{ext}$. 
We see that $\mathcal{C}_P^L$ is not rotation invariant, and that $\mathcal{C}_N^L$ and $\mathcal{C}_D^L$ are neither rotation nor translation invariant.
Moreover, since $\cov(\bfx, \bfy)>0$ for all $\bfx, \bfy \in \mathbb{R}^d$, we observe that
	\begin{equation}
		\cov_{P}^L(\bfx,\bfy) > \cov(\bfx,\bfy) \quad \text{ and }\quad  \cov_{N}^L(\bfx,\bfy) > \cov(\bfx,\bfy) \label{eq:pos_cov}.
	\end{equation}
That is, the folded covariance with either periodic or Neumann boundary conditions always overestimates the Mat\'ern covariance. 
This is not necessarily true for the Dirichlet case, since in \eqref{eq:dirsum} $\prod_{j=1}^d\varepsilon_j\in\left\{-1,1\right\}$.
\end{remark}

\section{Error estimate for the covariance on the bounded domain}\label{sect:mainthm}
A direct consequence of \cref{thm:folded} is the following result.
In what follows we write $\Matern(x)$, $x\ge 0$, to denote $\sigma^2 \mtn{\nu}(\kappa x)$ for convenience.
\begin{corollary}\label{cor:errsum}
For the covariances $\mathcal{C}^L_{\ast}$ for $\ast\in\left\{D,N,P\right\}$ the following error estimate holds
\begin{align}
	\lvert\cov_{\ast}^{L}(\bfx,\bfy) - \Matern(\bfx,\bfy)\rvert  &\leq  (2^d-1)\,\Matern(\delta) + 2^d\sum_{\bfk\in\mathbb{N}_0^d\setminus \left\{0\right\}}\Matern(L\lVert\bfk\rVert_2),\label{eq:errper}
	\end{align}
for all $\bfx,\bfy\in D$.
\end{corollary}
\begin{proof}
		
        We start with the estimate for periodic boundary conditions. 
	Without loss of generality, we \eu{assume that} all components of the vector~$(\bfx-\bfy)$ are positive. Indeed, if some components of $(\bfx-\bfy)$, with indices $j_1,\ldots,j_m$, $m\leq d$, are negative, then the right-hand side of \eqref{eq:cov_P} \eu{reads} 
	$$\sum_{\bfk\in\bbZ^d}\cov\left(\bfx+L\bfk,\bfy\right)=\sum_{\widetilde{\bfk}\in\bbZ^d}\cov\left(\widetilde{\bfx}+L\widetilde{\bfk},\widetilde{\bfy}\right),$$ where $\widetilde{x}_{i}=-x_{i}$, $\widetilde{y}_{i}=-y_{i}$, $\widetilde{k}_{i}=-k_{i}$ for $i=j_1,\ldots,j_m$, and $\widetilde{x}_{i}=x_{}$, $\widetilde{y}_{i}=y_{i}$, $\widetilde{k}_{i}=k_{i}$ otherwise.
	We denote by $\ev$ the vector with all components equal to~$1$. 
	Then, for any $\bfk\in\mathbb{Z}^d$, there exist $\bfk'\in\bbN_0^d$ and $\boldsymbol{\varepsilon}\in B^d$ such that
	\begin{equation}\label{eq:renumbering}
		\bfk = \bfeps.\left(\bfk'+\frac{1}{2}(\ev-\bfeps)\right).
	\end{equation}
	Hence, using \eqref{eq:cov_P},
	\begin{equation}
		\cov_P^{L}(\bfx,\bfy)
		= \sum_{\bfk\in\bbZ^d}\Matern(\bfx + \bfk L, \bfy)
		= \sum_{\bfk \in\bbN_0^d}\sum_{\bfeps\in B^d}\Matern(\abs{\bfk L + (\ev-\bfeps)\frac{L}{2} + \bfeps.(\bfx-\bfy)}).
		\label{eq:covP_inProof}
	\end{equation}
	Note that $\left((\ev-\bfeps)\frac{L}{2} + \bfeps.(\bfx-\bfy)\right) \in [0,L]^d$ for all $\bfx,\bfy\in D$.
	Then, due to the monotonicity of the Mat\'ern function, for all~$\bfeps\in B^d$ we have the bound
	\begin{equation}
		\sum_{\bfk \in\bbN_0^d}\Matern(\abs{\bfk L + (\ev-\bfeps)\frac{L}{2} + \bfeps.(\bfx-\bfy)})
		\leq 
		\sum_{\bfk \in\bbN_0^d\setminus\left\{\boldsymbol{0}\right\}}\Matern(\abs{\bfk } L)
		+ \Matern(\abs{(\ev-\bfeps)\frac{L}{2} + \bfeps.(\bfx-\bfy)}).
		\label{eq:covP_ineq}
	\end{equation}
	Moreover, if $\bfeps\ne\ev$, $\abs{(\ev-\bfeps)\frac{L}{2} + \bfeps.(\bfx-\bfy)} \ge \delta$.
	Thus, combining \eqref{eq:covP_inProof} and \eqref{eq:covP_ineq} we obtain
	\begin{equation}\label{1}
		\cov_P^{L}(\bfx,\bfy)
		\leq \Matern(\bfx,\bfy) + (2^d-1)\,\Matern(\delta) + 2^d\sum_{\bfk \in\bbN_0^d\setminus\left\{\boldsymbol{0}\right\}}\Matern(\abs{\bfk } L).
	\end{equation}
        Subtracting $\mathcal{C}(\bfx,\bfy)$ on both sides of \eqref{1}, and recalling \eqref{eq:pos_cov}, the desired error bound for periodic boundary conditions follows.

	For Neumann boundary conditions, we use again \eqref{eq:renumbering} to obtain
	\begin{align}
	\mathcal{C}_{N}^{L}(\bfx,\bfy)
	&= \sum_{\bfeps\in B^d}\sum_{\bfk\in\bbZ^d} \Matern\left(\bfx + 2L\bfk,\bfeps.\bfy\right) \nonumber\\
	&= \sum_{\bfeps_1,\bfeps_2\in B^d}\sum_{\bfk\in\bbN_0^d}
	\Matern\left(\abs{2L\bfk + (\ev-\bfeps_{1})L + \bfeps_1.(\bfx-\bfeps_2.\bfy)}\right). \label{eq:covN_inProof}
	\end{align}
	Note that the sum over $\bfk\in\bbZ^d$ is a sum over all cells $2L\bfk + [0,2L]^d$ (cf. \Cref{fig:folding}). For a fixed~$\bfeps_2$, let us consider $\bfz=\bfx-\bfeps_{2}.\bfy$ and address the summation over $\bfeps_1\in B^d$. We want to show that, in every cell $2L\bfk + [0,2L]^d$, $\bfk\in\bbZ^d$, the point $(\ev-\bfeps_{1})L + \bfeps_1.(\bfx-\bfeps_2.\bfy)$ belongs to  a different subcell of size $L$ as $\bfeps_1$ varies in $B^d$. 
	Indeed, for all~$\bfx,\bfy\in D$, we have $\bfz\in [0,2L]^d$. Hence, for all~$\bfeps_{1}$, it holds $\left((\ev-\bfeps_{1})L + \bfeps_{1}.\bfz\right) \in [0,2L]^d$. Moreover, observe that, for any component $\varepsilon_{1,j}\in B=\left\{-1,1\right\}$,
	\begin{align*}
		&&\left((1-\varepsilon_{1,j})L + \varepsilon_{1,j}z_j\right) &\in [0,L]
		& \Leftrightarrow &&
		\left((1+\varepsilon_{1,j})L - \varepsilon_{1,j}z_j\right) &\in [L,2L]. &&
	\end{align*}
	Therefore, summing over $\bfeps_1\in B^d$, we cover all $2^d$ subcells of size $L$ inside $2L\bfk + [0,2L]^d$.
	Using this fact and the monotonicity of the Mat\'ern function, for all $\bfeps_2\in B^d$ we obtain the estimate
	\begin{equation}\label{eq:covN_ineq}
		\sum_{\bfk\in\bbN_0^d}\sum_{\bfeps_1\in B^d}
		\Matern\left(\abs{2L\bfk + (\ev-\bfeps_{1})L + \bfeps_1.\bfz}\right)
		\le 
		\sum_{\bfk\in\bbN_0^d\setminus\left\{\boldsymbol{0}\right\}}
		\Matern\left(\abs{\bfk}L\right)
		+ \Matern\left(\abs{\min(\bfz,2L\ev - \bfz)}\right),
	\end{equation}
	where the minimum is applied element-wise.
	Observe also that, for $\bfx,\bfy\in D$, it holds
	\begin{equation*}
		\abs{\min(\bfx-\bfeps_{2}.\bfy,2L\ev - (\bfx-\bfeps_{2}.\bfy))} \ge \delta, \quad\text{if }\bfeps_2\ne\ev.
	\end{equation*}
	Thus, from \eqref{eq:covN_inProof} and \eqref{eq:covN_ineq} we obtain
	\begin{equation*}
		\cov_{N}^{L}(\bfx,\bfy)
		\le
		\Matern(\bfx,\bfy)
		+ (2^d-1)\,\Matern\left(\delta\right)
		+ 2^d\sum_{\bfk\in\bbN_0^d\setminus\left\{\boldsymbol{0}\right\}}\Matern\left(\abs{\bfk}L\right).
	\end{equation*}
	This and \eqref{eq:pos_cov} yields the desired error estimate for Neumann boundary conditions.
	
	Finally, let us consider Dirichlet boundary conditions.
	Using once again the renumbering in \eqref{eq:renumbering}, we write
	\begin{align}
		\cov_{D}^{L}(\bfx,\bfy)
		&= \sum_{\bfeps\in B^d}\left(\prod_{j=1}^{d}\bfeps_{j}\right)\sum_{\bfk\in\bbZ^d} \Matern\left(\bfx + 2L\bfk,\bfeps.\bfy\right) \nonumber\\
		&= \sum_{\bfeps_1,\bfeps_2\in B^d}\left(\prod_{j=1}^{d}\bfeps_{2,j}\right)\sum_{\bfk\in\bbN_0^d}
		\Matern\left(\abs{2L\bfk + (\ev-\bfeps_{1})L + \bfeps_1.(\bfx-\bfeps_2.\bfy)}\right). \label{eq:covD_inProof}
	\end{align}
	The inequality \eqref{eq:covN_ineq} still holds.
	Using this and the fact that $\cov(x)>0$ to bound the terms with signs $\prod_{j=1}^{d}\bfeps_{2,j} = 1$ and $\prod_{j=1}^{d}\bfeps_{2,j} = -1$ respectively, and then vice versa, we obtain the upper and lower bounds
	\begin{align*}
	\mathcal{C}_{D}^{L}(\bfx,\bfy)
	&\le
	\Matern(\bfx,\bfy)
	+ (2^{d-1}-1)\,\Matern\left(\delta\right)
	+ 2^{d-1}\sum_{\bfk\in\bbN_0^d\setminus\left\{\boldsymbol{0}\right\}}\Matern\left(\abs{\bfk}L\right), \\
	\cov_{D}^{L}(\bfx,\bfy)
	&\ge
	\Matern(\bfx,\bfy)
	- 2^{d-1}\,\Matern\left(\delta\right)
	- 2^{d-1}\sum_{\bfk\in\bbN_0^d\setminus\left\{\boldsymbol{0}\right\}}\Matern\left(\abs{\bfk}L\right).
	\end{align*}
	Hence,
	\begin{equation}\label{eq:better}
		\lvert\cov_{D}^{L}(\bfx,\bfy)-\Matern(\bfx,\bfy)\rvert
		\le
		2^{d-1}\,\left(\Matern\left(\delta\right)
		+\sum_{\bfk\in\bbN_0^d\setminus\left\{\boldsymbol{0}\right\}}\Matern\left(\abs{\bfk}L\right)\right).
	\end{equation}
	Observe that this bound is smaller than the bound in the statement of \eqref{cor:errsum}.
\end{proof}

In \eqref{eq:errper} we can distinguish two error contributions, one that depends on $\delta$ and one that depends on the size $L$ of $D_{ext}$. 
Due to the exponential decay of the Mat\'ern covariance, we expect the second contribution to be negligible for $L$ sufficiently large. 
We prove this next.

\begin{theorem}[Main result]\label{thm:mainres}
For $\bfx,\bfy\in D$, the error in the covariance when using homogeneous Neumann, homogeneous Dirichlet or periodic boundary conditions on $D_{ext}$ is bounded by
		\begin{equation}\label{eq:Error1}
		\lVert\cov_{\ast}^{L}(\bfx,\bfy) - \Matern(\bfx,\bfy)\rVert_{\infty} 
		\leq A\cdot \sigma^2\mtn{\nu}(\kappa\delta),
		\end{equation}
		$\ast\in\left\{D,N,P\right\}$. The constant~$A$ is given by
		\begin{equation}\label{eq:A}
		A = (2^d-1)\cdot\left(1 + \frac{2^d d!\cdot f(\ell)}{(1-f(\ell))^d}\right),\quad \ell=L-\delta,
		\end{equation}
		where the function~$f(x)$ is defined for any $x>0$ as
		\begin{equation}\label{eq:ffunc}
		f(x) = \mtn{\max(\nu,1/2)}(\kappa x)
		\end{equation}
(with $\kappa$ depending on $\nu$ as in \eqref{eq:matcov}).
\end{theorem}

\begin{proof}
In this proof, we use the auxiliary result in \cref{lem:covbound} in the Appendix.
Recall that for convenience we use the notation $\Matern(x)=\sigma^2 \mtn{\nu}(\kappa x)$, $x\ge 0$.
According to \cref{cor:errsum}, it is sufficient to bound the last summand on the right-hand side of \eqref{eq:errper}. 
To this end, we observe that
\begin{align}
	\sum_{\bfk\in\bbN_0^d\setminus\left\{\boldsymbol{0}\right\}}\Matern(\abs{\bfk} L)&\leq \sum_{\bfk\in\bbN_0^d\setminus\left\{\boldsymbol{0}\right\}}\Matern(\lVert \bfk\rVert_\infty L)  \leq d\sum_{k=1}^{\infty}(k+1)^{d-1}\cdot \Matern(kL)  \nonumber\\
	& \leq 2^{d-1} d \sum_{k=1}^{\infty}k^{d-1}\cdot \Matern(kL),   \label{eq:Sum}
\end{align}
where for the first inequality we have used the monotonicity of the Mat\'ern kernel, and in the second one the fact that, for every $k\in\bbN$, $\sharp\left\{\bfk\in\bbN_0^d: \lVert\bfk\rVert_{\infty}=k\right\}\le d(k+1)^{d-1}$.
Using the function \eqref{eq:ffunc}, we have from \cref{lem:covbound} that $\Matern(kL)\leq\Matern(L)f(L)^{k-1}$, for every $k\in\bbN$, therefore
\begin{equation}\label{eq:bound1}
	\sum_{\bfk\in\bbN_0^d\setminus\left\{\boldsymbol{0}\right\}}\Matern(\abs{\bfk} L)\leq \Matern(L) \cdot \left( 2^{d-1} d\sum_{k=1}^{\infty}k^{d-1}\cdot f(L)^{k-1} \right).
\end{equation}
We note that $0<f(L)<1$. 
The sum on the right-hand side of \eqref{eq:bound1} involves the polylogarithm $\text{Li}_s(z)=\sum_{k=1}^{\infty}k^{-s}z^k$, for $s,z \in\bbC$, and it can be bounded as
	\begin{align}
	\sum_{k=1}^{\infty}k^{d-1}\cdot f(L)^{k-1} \nonumber
	&= \frac{1}{f(L)}\Li{-(d-1)}(f(L))= \frac{1}{(1-f(L))^d}\sum_{k=0}^{d-2}\EuNum{d-1}{k}f(L)^{d-2-k} \nonumber\\
	&< \frac{(d-1)!}{(1-f(L))^d}, \label{eq:polylog}
	\end{align}
where $\EuNum{n}{k}$ for $n, k\in\bbN$ are the Eulerian numbers, and we have used that $\sum_{k=0}^{n-1}\EuNum{n}{k} = n!$. Inserting \eqref{eq:polylog} in \eqref{eq:bound1} and using that \cref{lem:covbound} gives $\Matern(L)=\Matern(\ell+\delta)\leq \Matern(\delta)f(\ell)$, we obtain
	\begin{align*}
	\sum_{\bfk\in\bbN_0^d\setminus\left\{\boldsymbol{0}\right\}}\Matern(\abs{\bfk} L) 
	&\leq  \frac{2^{d-1}d!}{(1-f(L))^d}\cdot \Matern(L)\leq \frac{2^{d-1}d!\cdot f(\ell)}{(1-f(\ell))^d}\cdot \Matern(\delta).
	\end{align*}
The above inequality together with \cref{cor:errsum} leads then to the desired result:
	\begin{align*}
	\lVert \cov_{*}^{L}(\bfx,\bfy) - \Matern(\bfx,\bfy)\rVert_\infty
	&\le (2^d-1)\,\Matern(\delta) + 2^d\sum_{\bfk\in\bbN_0^d\setminus\left\{\boldsymbol{0}\right\}}\Matern(\abs{\bfk} L) \\
	&\le (2^d-1)\,\Matern(\delta) + 2^d\cdot \frac{2^{d-1}d!\cdot f(\ell)}{(1-f(\ell))^d}\cdot \Matern(\delta) \\
	&\le (2^d-1)\left(1 + \frac{2^d d!\cdot f(\ell)}{(1-f(\ell))^d}\right) \cdot \sigma^2\mtn{\nu}(\kappa\delta). \\
	\end{align*}
\end{proof}

\begin{remark}
The estimate in \eqref{eq:Error1} shows that the error behaves as $\mtn{\nu}(\kappa\delta)$. 
In particular, due to the definition of $\kappa$ (see \eqref{eq:matcov}), the magnitude of the error depends on the ratio $\eu{\delta / \rho}$, and, asymptotically, the error decreases exponentially as this ratio increases. 
This explains the observations in \cite{LRL} and \cite{Roininen14}, namely, that the error in the covariance is negligible if the distance of $\partial D_{ext}$ from $\partial D$ is greater than the correlation length $\rho$, i.e. $\eu{\delta / \rho} > 2$.
\end{remark}

\begin{remark}
Using a scaling argument, it is easy to see that the error bound in \cref{thm:mainres} also holds for rectangular bounding boxes. 
We consider, for instance, periodic boundary conditions on the hyperrectangle $\eu{\bigtimes_{i=1}^d} [0,L_i]$, for $L_i>0$, $i=1,\ldots,d$. 
Defining $\widetilde{\bfk}:=\left(\frac{k_1}{L_1},\ldots,\frac{k_d}{L_d}\right)$ for every $\bfk\in\bbN^d$, the eigenpairs are
\begin{equation*}
w_{\bfk}= \frac{1}{\left(L_1\ldots L_d\right)^{\frac{1}{2}}}e^{2\pi\img\tilde{\bfk}\cdot\bfx},\quad 
\lambda_{\bfk}=1+\left(\frac{2\pi}{\kappa}\right)^2 \lVert\tilde{\bfk}\rVert_2^2,
\end{equation*}
$\bfk\in\bbN_0^d$, cf. \eqref{eq:eigper}. 
The expression for the folded covariance \eqref{eq:cov_P} is, in this case,
\begin{equation*}
\mathcal{C}_{P}^{\bfL}(\bfx,\bfy) = \sum_{\bfk\in\mathbb{Z}^d}\mathcal{C}(\bfx+\bfL.\bfk,\bfy),
\end{equation*}
with $\bfL:=\left(L_1,\ldots,L_d\right)$ (\eu{recall} that $\bfL.\bfk$ \eu{denotes} the component-wise product). 
Then,  \eu{the estimate in} \eqref{eq:errper} is modified as
\begin{equation*}
	\lvert\cov_{\ast}^{\bfL}(\bfx,\bfy) - \Matern(\bfx,\bfy)\rvert  \leq  (2^d-1)\,\Matern(\delta) + 2^d\sum_{\bfk\in\mathbb{N}_0^d\setminus \left\{0\right\}}\Matern(\lVert\bfL.\bfk\rVert_2),
\end{equation*}
and the result of \cref{thm:mainres} holds by \eu{letting} $\ell=\min_{i=1,\ldots,d}L_i-\delta$.
\end{remark}

\section{Further analyses}\label{sect:ext}
\eu{In this section we discuss the extension of the main result in \Cref{thm:mainres} to \textit{anisotropic} Mat\'ern kernels.
In addition we comment on the use of Robin boundary conditions.}

\newcommand{\bfTheta}{\mathbf{\Theta}}

\subsection{Anisotropic Mat\'ern covariance}
	One can consider a more general form of Mat\'ern covariance (see e.g. \cite{Stein2005}):
	\begin{equation}\label{eq:matcov_anis}
	\cov(\bfx,\bfy) 
	= 
	\sigma^2 \mathcal{M}_\nu(\sqrt{2\nu}\,\Vert \bfx-\bfy\Vert_{\bfTheta^{-1}}),\quad \Vert \bfx-\bfy\Vert_{\bfTheta^{-1}} = \sqrt{(\bfx-\bfy)\cdot\bfTheta^{-1}\cdot(\bfx-\bfy)},
	\end{equation}
	with the metric given by the constant, symmetric positive definite matrix~$\bfTheta = \mathbf{R}\mathbf{D}^2\mathbf{R}^{-1}$, where $\mathbf{D}$ is a diagonal matrix with entries $\rho_j>0$, $j=1,\dots,d$, and $\mathbf{R}$ is an orthogonal matrix.	
	Introducing the coordinate transformation $\widetilde{\bfx} = (\mathbf{R}\mathbf{D})^{-1}\bfx$ (scaling and rotation), we obtain the form \eqref{eq:matcov} with $\rho=1$.
	Hence, \eqref{eq:matcov_anis} writes through its Fourier integral:
	\begin{align*}
	\cov(\bfx,\bfy) 
	&= \frac{\eta^2}{(2\pi)^d}\int_{\bbR^d}  \left(1 + \frac{1}{2\nu}\lVert\bfxi\rVert_2^2 \right)^{-\alpha}\, e^{\img \bfxi \cdot(\mathbf{R}\mathbf{D})^{-1}(\bfx-\bfy)}  \dd\bfxi \\
	&= \frac{\eta^2}{(2\pi)^d}\int_{\bbR^d}  \left(1 + \frac{1}{2\nu}\lVert\widetilde{\bfxi}\rVert_\bfTheta^2 \right)^{-\alpha}\, e^{\img \widetilde{\bfxi} \cdot(\bfx-\bfy)}  (\det\mathbf{D})\dd\widetilde{\bfxi}, 
	& \eta^2 = \sigma^2\frac{\left(2\pi\right)^{d/2}\Gamma(\nu+d/2)}{\nu^{d/2}\,\Gamma(\nu)},
	\end{align*}
	where we used the change of variables $\widetilde{\bfxi} = \mathbf{R}\mathbf{D}^{-1} \bfxi$, and the $\bfTheta$-norm $\Vert \widetilde{\bfxi}\Vert_{\bfTheta}^2 = \widetilde{\bfxi} \cdot \bfTheta\widetilde{\bfxi}$.
	Thus, the SPDE for \eqref{eq:matcov_anis} is
	\begin{equation}\label{eq:spde_anis}
	\left(\mathcal{I}-\frac{1}{2\nu}\,\nabla\cdot\left(\mathbf{\Theta}\nabla\right)\right)^{\frac{\alpha}{2}} u(\omega,\bfx) 
	= 
	\widehat{\eta}\dot{W}(\omega,\bfx),\quad \bfx\in\mathbb{R}^d,\quad \text{for }\mathbb{P}\text{-a.e. } \omega\in\Omega,
	\end{equation}
	with
	\begin{equation*}
	\widehat{\eta}^2 = \sigma^2\sqrt{\det\bfTheta}\cdot\frac{\left(2\pi\right)^{d/2}\Gamma(\nu+d/2)}{\nu^{d/2}\,\Gamma(\nu)}.
	\end{equation*}
	The eigenfunctions of the operator~$\left(\mathcal{I}-\frac{1}{2\nu}\,\nabla\cdot\left(\mathbf{\Theta}\nabla\right)\right)$ on~$D_{ext}$ with homogeneous Neumann, Dirichlet or periodic boundary conditions are still given by \eqref{eq:eigneu}, 
	(\ref{eq:eigdir}) and (\ref{eq:eigper}), respectively.
	However, the eigenvalues now involve the $\bfTheta$-norm, i.e. they read
	\begin{equation*}
		\lambda_{\bfk}=1 + \frac{1}{2\nu}\,\left(\frac{\pi}{L}\right)^2\lVert\bfk\rVert_\bfTheta^2.
	\end{equation*}
	for Dirichlet and Neumann boundary conditions, and analogously for periodic boundary conditions.
	The proof of \cref{thm:mainres} can be reproduced, given $$\lVert\bfx\rVert_{\bfTheta^{-1}} \ge \rho_{max}^{-1}\lVert\bfx\rVert_2 \ge {\rho_{max}^{-1}}\lVert\bfx\rVert_\infty,$$ where $\rho_{max} = \max\limits_{j=1,\dots,d}\rho_j$, that yields the error bound
	\begin{equation}
	\lVert\cov_{\ast}^{L}(\bfx,\bfy) - \Matern(\bfx,\bfy)\rVert_{\infty} 
	\leq A\cdot \sigma^2\mtn{\nu}\left(\frac{\sqrt{2\nu}}{\rho_{max}}\delta\right),
	\end{equation}
	where $\ast\in\left\{D,N,P\right\}$, $A = (2^d-1)\cdot\left(1 + \frac{2^d d!\cdot f(\ell)}{(1-f(\ell))^d}\right)$, $\ell=L-\delta$, and 	$f(x) = \mtn{\max(\nu,1/2)}\left(\frac{\sqrt{2\nu}}{\rho_{max}} x\right)$.
	\eu{SPDEs with anisotropic Laplacian as in \eqref{eq:spde_anis} are studied in \cite{Fuglstad2015}}.

\subsection{Robin boundary conditions}\label{subsect:robin}
We have mentioned in the introduction that another possibility is to use Robin boundary conditions $\nabla u_L \cdot \bfn +\beta u_L =0$ on $D_{ext}$, $\beta>0$. 
In order to get boundary conditions minimizing the boundary effects, the coefficient $\beta$ needs to be tuned. 
In \cite{Daon2018} it is \eu{shown} that $\beta=\kappa$ is the exact boundary condition for $d=1$ and $\nu=0.5$. This is because, in this case, the Mat\'ern kernel corresponds to the exponential function, satisfying these Robin conditions. 

For $d=1$ and $\nu \neq 0.5$,  and for $d>1$ and any $\nu>0$, $\beta=\kappa$ does not provide exact boundary conditions. 
However, such choice for the Robin coefficient can be motivated by analogies with absorbing boundary conditions for the Helmholtz equation. Indeed, for every $\bfx\in D$ fixed, and denoting $r=\lVert\bfy-\bfx\rVert_2$ for every $\bfy\in\bbR^d$, the Mat\'ern covariance satisfies 
$$\left(\kappa r\right)^{-\nu+\frac{1}{2}}\left(\dfrac{\partial}{\partial r} + \kappa\right)  \mtn{\nu}\left(\kappa r\right)\sim  e^{-\kappa r}, \quad\text{for }r\rightarrow\infty,$$
with a constant depending on $\nu$.
\eu{In particular, for every $\nu>0$ it holds}
\begin{equation}\label{eq:rad}
\lim_{r\rightarrow\infty} \left(\kappa r\right)^{-\nu+\frac{1}{2}}\left(\dfrac{\partial}{\partial r} + \kappa\right) \mtn{\nu}\left(\kappa r\right) =0.
\end{equation}
This can be seen by computing $$\left(\dfrac{\partial}{\partial r} + \kappa\right) \mtn{\nu}\left(\kappa r\right)=\kappa\left(\kappa r\right)^{\nu}\left(\BesselK{\nu}(\kappa r)-\BesselK{\nu-1}(\kappa r)\right)$$ and using the asymptotic expansion for the modified Bessel functions \cite[Sect. 7.23]{watson1995treatise}. 

Equation \eqref{eq:rad} can be seen as the correspondent, in our case, of what for the Helmholtz equation is the Sommerfeld radiation condition \cite[Sect. 2.1]{CK}. For the Helmholtz equation, first order absorbing boundary conditions, corresponding to Robin boundary conditions, are obtained by imposing the Sommerfeld condition at the boundary of the domain \cite{BabStro08}. If for \eqref{eq:spde} we impose the boundary conditions $\nabla u_L \cdot \bfn +\kappa u_L =0$ on $D_{ext}$, then the associated covariance is the Green's function of $\left(\mathcal{I}+\kappa^{-2}\Delta\right)^{\alpha}$ associated to these Robin boundary conditions, and this means imposing \eqref{eq:rad} at finite distance, 
on $\partial D_{ext}$. In \Cref{sect:numexp}, we show the performance of the Robin boundary condition $\nabla u_L \cdot \bfn +\kappa u_L =0$ on $D_{ext}$ for $d=1,2$ and different values of $\nu$.


\section{Numerical results}\label{sect:numexp}
We present here numerical experiments to verify the error bound in \cref{thm:mainres}, for $d=1$ and $d=2$. 
Also, we compare Dirichlet, Neumann and periodic boundary conditions with the Robin boundary conditions $\nabla u_L\cdot\bfn + \kappa u_L=0$ on $\partial D_{ext}$ as from \Cref{subsect:robin}.
We consider normalized covariances, that is $\sigma^2=1$. 
For all boundary conditions, we use the analytic expression \eqref{eq:covexp} for the covariance, with \eqref{eq:eigdir}, \eqref{eq:eigneu} and \eqref{eq:eigper} for Dirichlet, Neumann and periodic boundary conditions, respectively. 
For Robin boundary conditions, we can also use \eqref{eq:covexp}, with eigenvalues and eigenfunctions given by equations (3.2)--(3.4) in \cite{Lapleigs13}.\footnote{In \cite{Lapleigs13}, equation (3.4) has a minor typo and the norm of eigenfunctions is $\lVert u_n^{(i)}(x)\rVert^2_{L^2((0,\ell_i))}=\left(\frac{\alpha_n^2+2h\ell_i+h^2\ell^2_i}{2h^2 \ell_i}\right).$}
The infinite sum in \eqref{eq:covexp} is truncated at $\lVert\bfk\rVert_{\infty}=\lceil \frac{L}{h} \rceil + 1$, with $h$ chosen sufficiently small to guarantee that the truncation error is negligible compared to the error in the covariance (ranging from $h=1\cdot 10^{-6}$ to $h=5\cdot 10^{-3}$). 
To evaluate the Mate\'rn covariance we use the formula in \eqref{eq:matcov}. 
The error \eqref{eq:linferr} has been computed taking the maximum over all pairs $(\bfx,\bfy)$ belonging to a discrete grid, specified below.

\subsection{Experiments for $d=1$}
\begin{figure}
	\centering
		\includegraphics[height=0.25\textwidth]{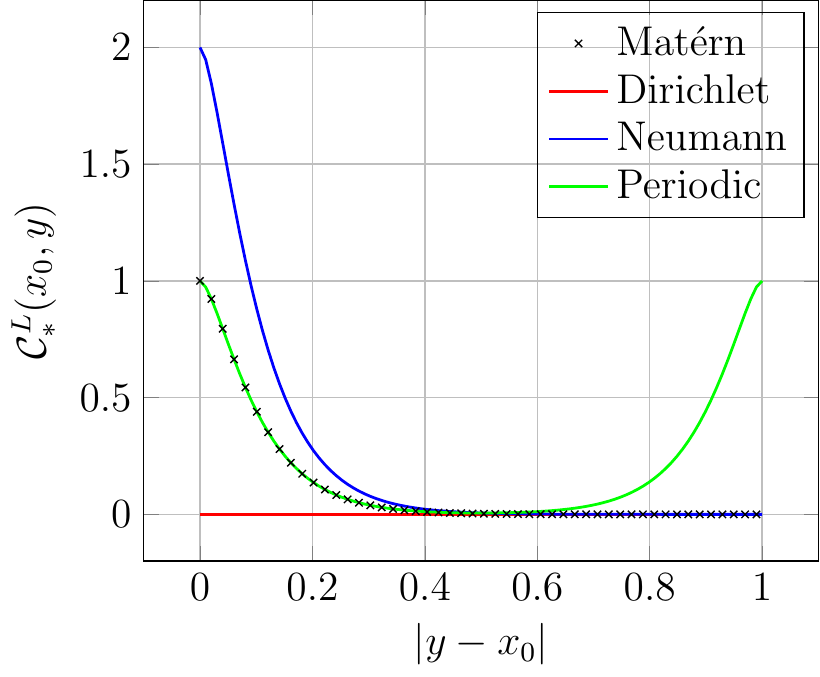}
	\hfill
	\includegraphics[height=0.25\textwidth]{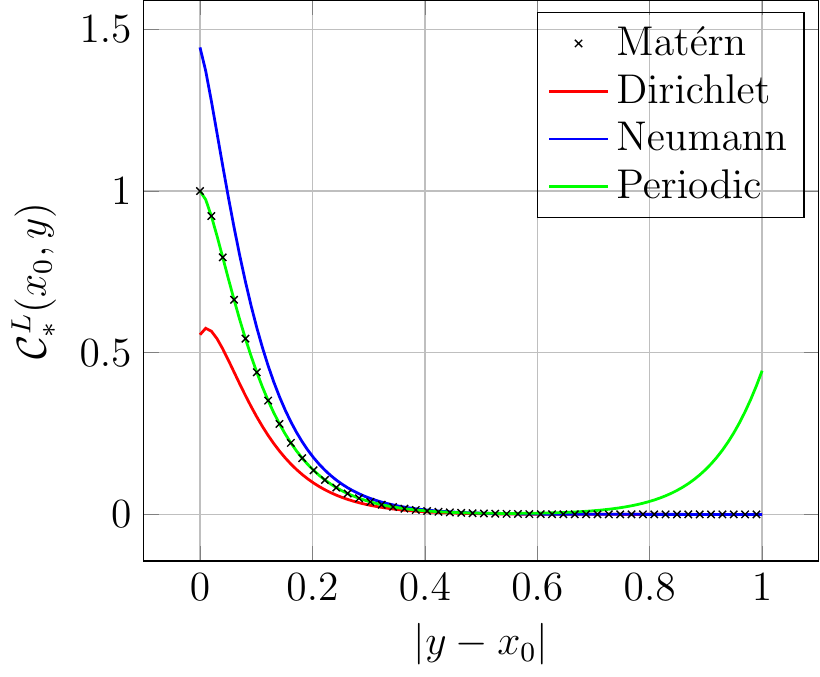}
	\hfill
	\includegraphics[height=0.25\textwidth]{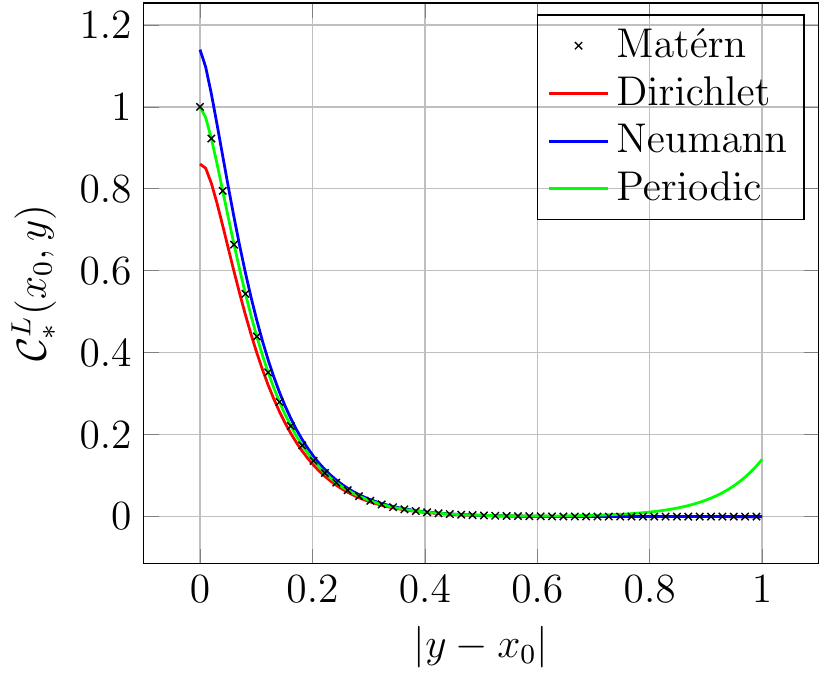}
	\caption{Case $d=1$, $\nu=1$ and $\rho=0.1$. Covariance functions $\cov^L_{\ast}(x_0,y)$, $\ast\in\left\{P,D,N\right\}$, in comparison with the Mat\'ern kernel. Here $x_0=\frac{\delta}{2}$ and $y\in[\frac{\delta}{2},1+\frac{\delta}{2}]$. Left: $\delta=0$. Center: $\delta=\rho$. Right: $\delta=2\rho$.}\label{fig:1dcovdifferentbc}
\end{figure}

\begin{figure}
	\centering
		\includegraphics[height=0.25\textwidth]{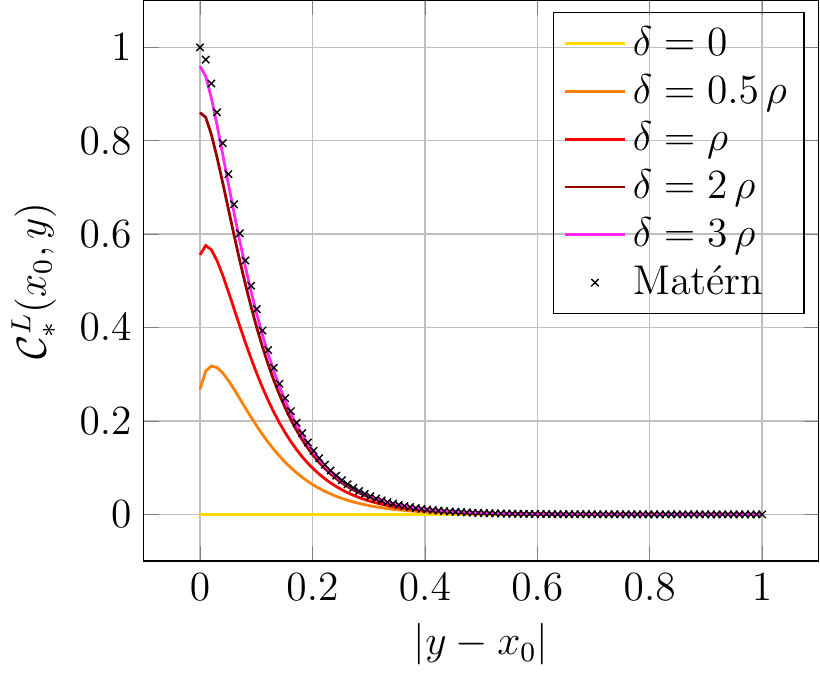}
	\hfill
	\includegraphics[height=0.25\textwidth]{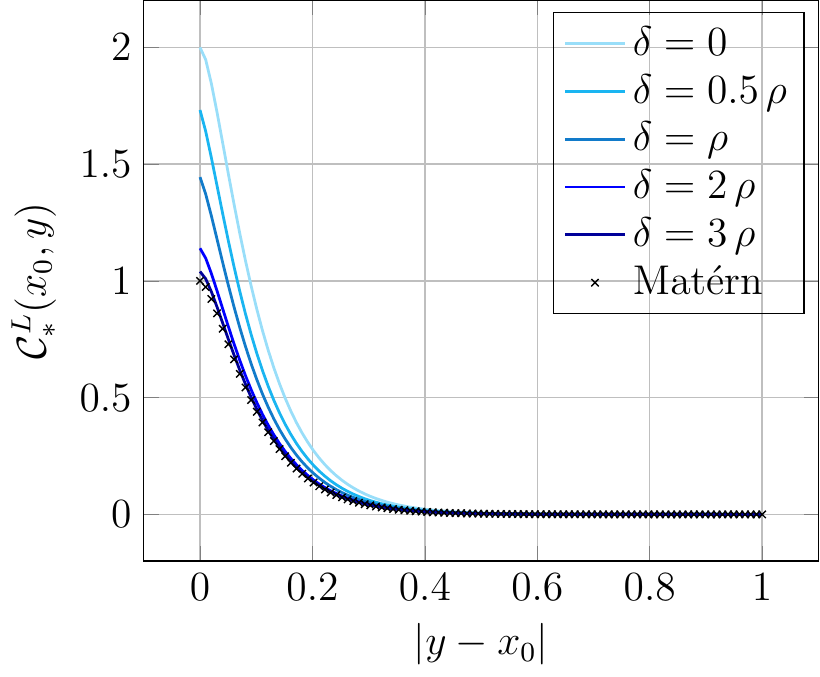}
	\hfill
	\includegraphics[height=0.25\textwidth]{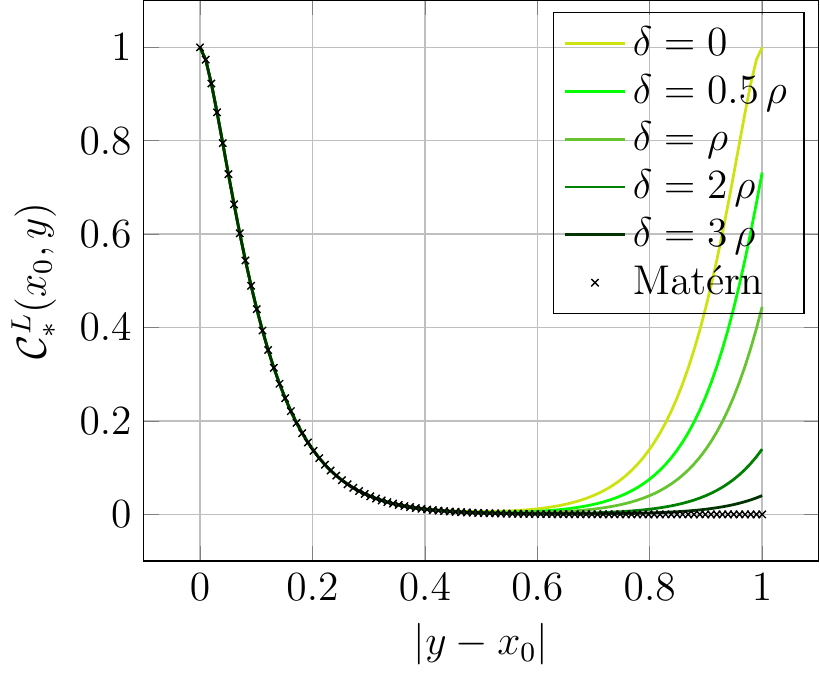}
	\caption{Case $d=1$, $\nu=1$ and $\rho=0.1$. Covariance functions $\cov^L_{\ast}(x_0,y)$, $\ast\in\left\{P,D,N\right\}$, in comparison with the Mat\'ern kernel, for different window sizes. Here $x_0=\frac{\delta}{2}$ and $y\in[\frac{\delta}{2},1+\frac{\delta}{2}]$. Left: Dirichlet b.c. ($\ast=D$). Center: Neumann b.c. ($\ast=N$). Right: periodic b.c. ($\ast=P$).}\label{fig:varydelta}
\end{figure}

For these experiments, we consider the domain $D=\left(\frac{\delta}{2},1+\frac{\delta}{2}\right)$, with $\frac{\delta}{2}\geq 0$ the window size.
\cref{fig:1dcovdifferentbc} shows  the effect of different boundary conditions on the covariance for $\nu=1$ and $\rho=0.1$. 
Although $\cov^L_{\ast}(x,y)$ for $\ast\in\left\{D,N\right\}$ depend on the position of $x,y\in D$ and not only on their distance, \cref{fig:1dcovdifferentbc} reports only the covariances $\cov^L_{\ast}(x_0,y)$, $\ast\in\left\{P,D,N\right\}$, and the Mat\'ern covariance for $x_0=\frac{\delta}{2}$ fixed and $y\in [\frac{\delta}{2},1+\frac{\delta}{2}]$. 
Apart from the expected result that, as $\delta$ increases, the approximate covariances become closer to the Mat\'ern covariance, we can also observe a different behavior for Dirichlet or Neumann and periodic boundary conditions: while the Dirichlet and Neumann boundary conditions affect mostly short range correlations, the periodic ones affect mostly long range correlations. 
This has to be expected from the nature of periodic boundary conditions. 
Indeed, we have seen in \eqref{eq:cov_N} and \eqref{eq:cov_D} that the covariances for Dirichlet and Neumann boundary conditions are the sum of covariances for periodic boundary conditions but with twice the period. 
Therefore, these boundary conditions introduce long range correlations, but for distances larger than the size of the domain of interest. 
How the approximate covariance approaches the Mat\'ern covariance as the window size increases can be observed, qualitatively, in \cref{fig:varydelta}, for the three types of boundary conditions. 

\begin{figure}
	\centering
	\begin{tabular}{ll}
		\includegraphics[height=0.25\textwidth]{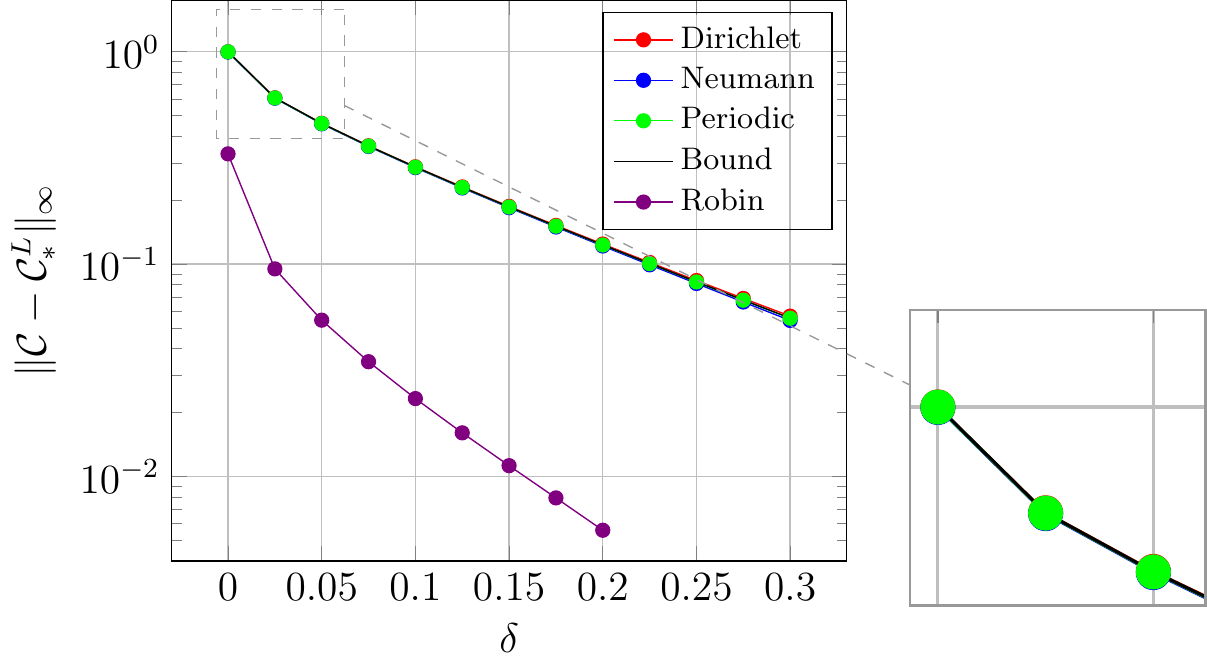} &	\includegraphics[height=0.25\textwidth]{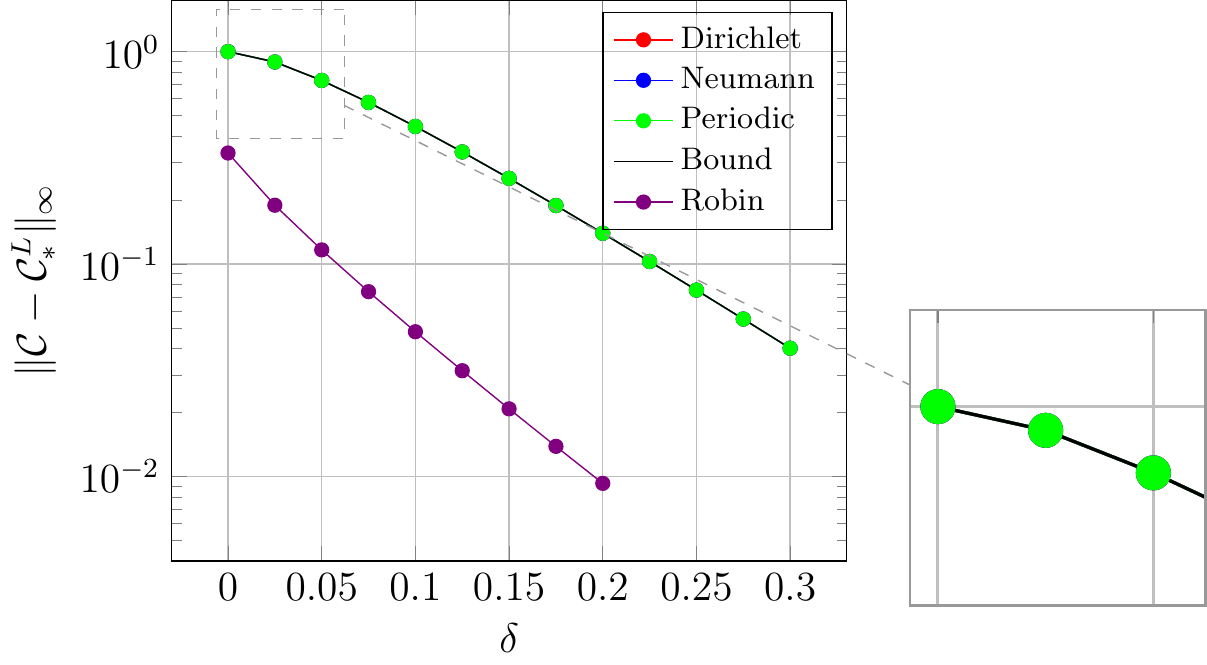}\\
	\\
	\includegraphics[height=0.25\textwidth]{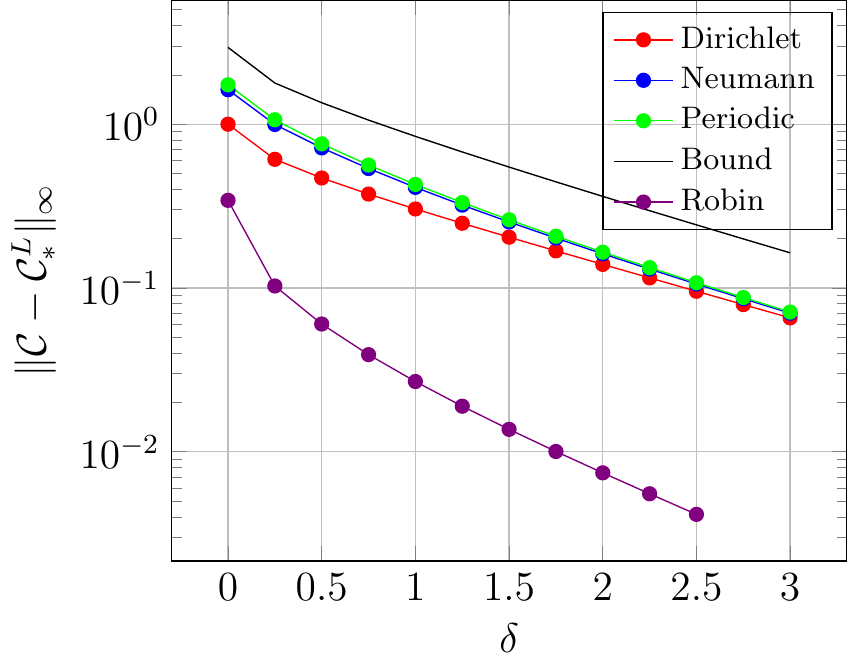} &	\includegraphics[height=0.25\textwidth]{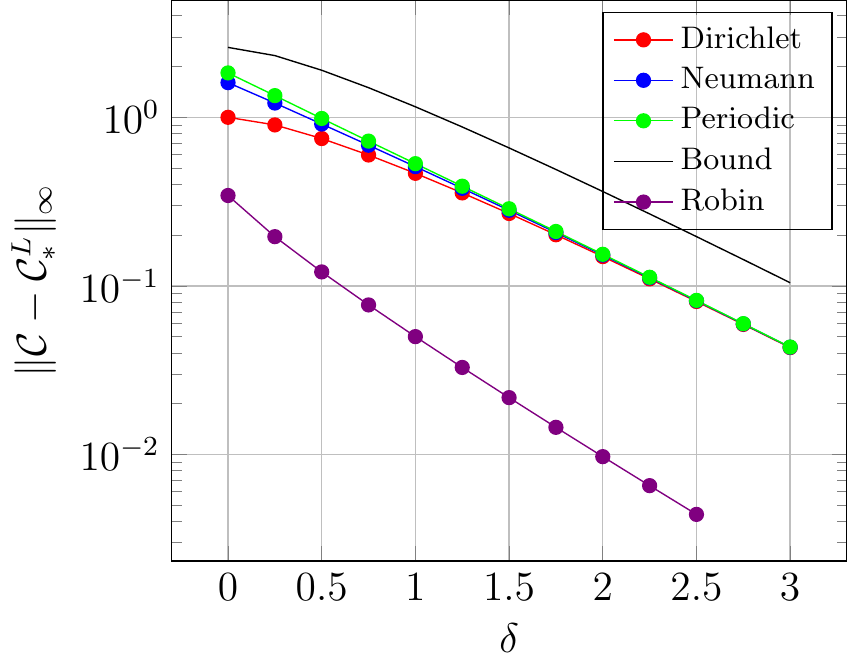}
	\end{tabular}
	\caption{Case $d=1$. Maximum norm error in the covariance kernel as function of the window size~$\delta$ for different boundary conditions. Left column: $\nu=0.25$. Right column: $\nu=1$. Top row: $\rho=0.1$. Bottow row: $\rho=1$.}\label{fig:bound1d}
\end{figure}

The bound stated in \cref{thm:mainres} is verified in \cref{fig:bound1d}. 
There, the error \eqref{eq:linferr} has been approximated by $\max_{x,y\in \mathcal{G}}\left|\Matern(x,y)-\cov^L_{\ast}(x,y)\right|$, for $\mathcal{G}$ a grid of $n$ equispaced points in $\overline{D}$ (the covariances are continuous functions, therefore the norm in \eqref{eq:linferr} coincides with the $C^0(\overline{D}\times\overline{D})$-norm). 
We have used $n=15$ for $\nu=1$ and $n=10$ for $\nu=0.25$, because for the latter case we need more terms in the spectral expansion for the covariance and it is therefore computationally more intense.
We can observe that the bound \eqref{eq:Error1} holds for Dirichlet, Neumann and periodic boundary conditions and both $\nu<0.5$ and $\nu\geq 0.5$. 
If the correlation length is significantly smaller than the domain size, see top row of \cref{fig:bound1d} for $\rho=0.1$, then the bound is very sharp and all boundary conditions, while behaving differently as observed in \cref{fig:1dcovdifferentbc}, produce the same error in the maximum norm. 
As the correlation length increases, the error bound \eqref{eq:Error1} is not as sharp.
Moreover, Dirichlet boundary conditions provide smaller errors than Neumann and periodic conditions in the preasymptotic regime. 
This is observed in the bottom row of \cref{fig:bound1d} for the extreme case of a correlation length equal to the size of the domain. 
Comparing with Robin boundary conditions, we can observe that these ones deliver lower errors in all cases considered. 
We have also tested that these conditions are exact for $\nu=0.5$, obtaining an error in the covariance stemming from the truncation of \eqref{eq:covexp} only.

\begin{figure}
	\centering
	\includegraphics[height=0.28\textwidth]{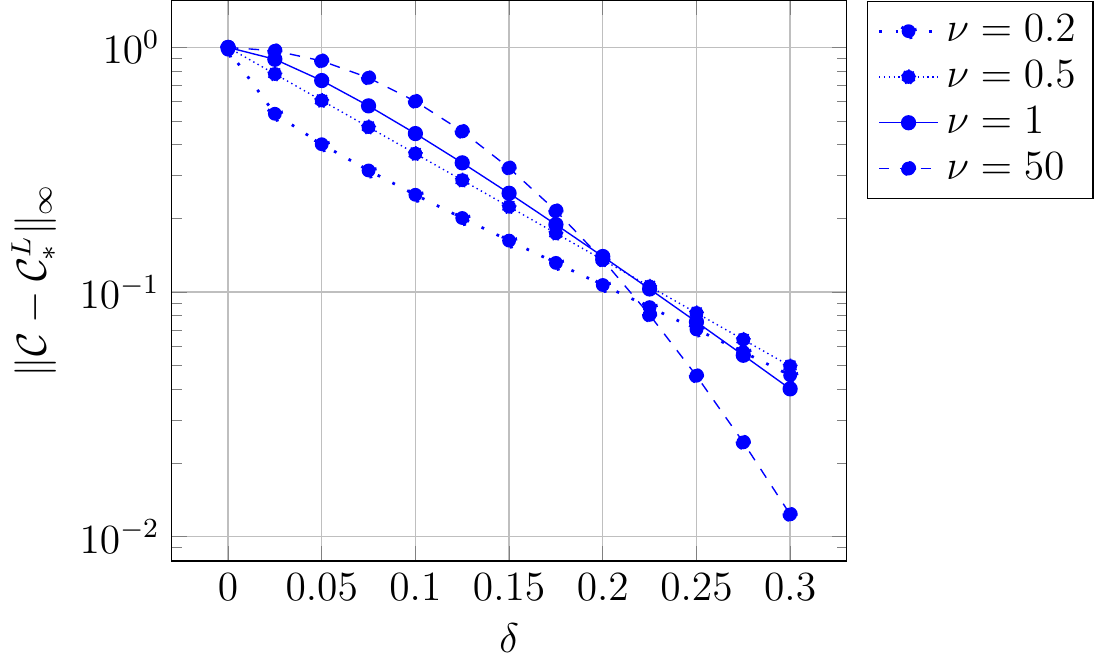}
	\hfill
	\includegraphics[height=0.28\textwidth]{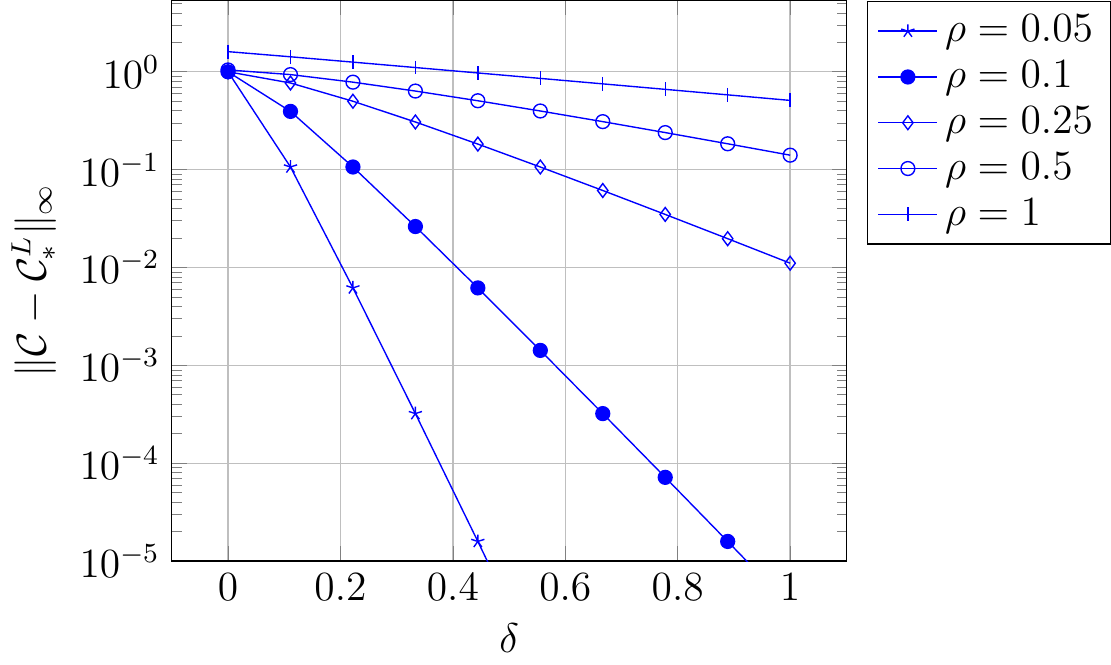}
	\caption{Case $d=1$, Neumann boundary conditions. Maximum norm error in the covariance kernel as function of the window size~$\delta$. Left: error for different values of~$\nu$ and fixed~$\rho=0.1$.
		Right: error for different values of~$\rho$ and fixed~$\nu=1$.}\label{fig:varynurho}
\end{figure}

Next, we fix the boundary conditions to the homogeneous Neumann ones and observe how the error behaves as $\nu$ or $\rho$ vary. 
The results for $\rho=0.1$ and $\nu$ varying are reported in the left plot of \cref{fig:varynurho}. 
For all values of $\nu$, the result of \cref{thm:mainres} is verified, and we see that the error behaves as $\mtn{\nu}(\sqrt{2\nu}\,\delta/\rho)$ (indeed, we have not drawn the bounds from \eqref{eq:Error1} as  they overlap with the error curves).
For $\nu=0.5$, we can observe a straight line in the semilogarithmic plot, reflecting the fact that for $\nu=0.5$ the Mat\'ern kernel corresponds to the exponential kernel $\mtn{\frac{1}{2}}(x)=e^{-\kappa x}$, $x\ge0$. 
For $\nu<\frac{1}{2}$, the logarithmic error curves are concave, and for $\nu>\frac{1}{2}$ they are convex.
Overall, for fixed $\rho$, the slopes of tangents of the error curves at $\delta=0$ decrease as $\nu$ increases, while, for $\delta\rightarrow\infty$, the slopes increase as $\nu$ increases, since the error behaves as $e^{-\kappa \delta}$.
 
If instead we fix $\nu=1$ and vary the correlation length $\rho$, we obtain the results in the right plot of \cref{fig:varynurho}. Since from \cref{fig:varydelta} it is clear that, for Neumann boundary conditions, the maximum error occurs when $x=y$, to compute the maximum norm we have used a grid $\mathcal{G}$ with $n=2$ points. 
In the right plot of \cref{fig:varynurho}, the case $\rho=1$ corresponds to the extreme case where the correlation length is equal to the size of the domain. 
Since asymptotically the error behaves as $e^{-\frac{\sqrt{2\nu}}{\rho}\delta}$, the error curves become steeper as $\rho$ increases.

\subsection{Experiments for $d=2$} 
\begin{figure}
	\centering
	\includegraphics[height=0.25\textwidth]{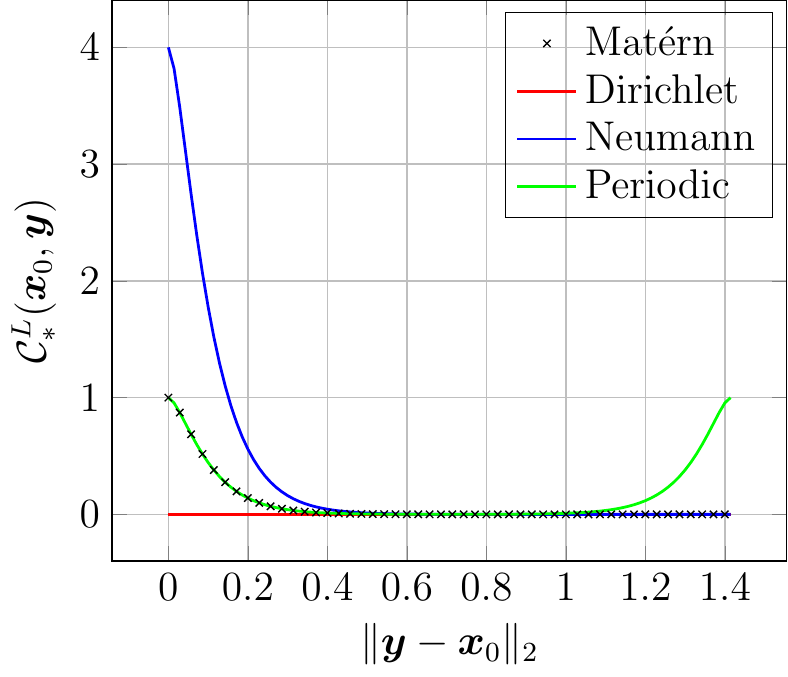}
	\hfill
	\includegraphics[height=0.25\textwidth]{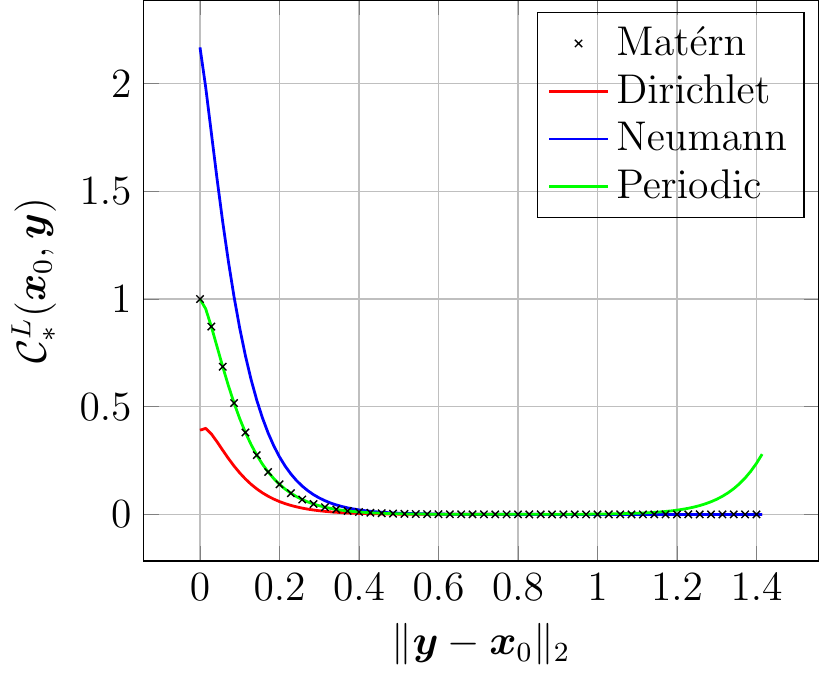}
	\hfill
	\includegraphics[height=0.25\textwidth]{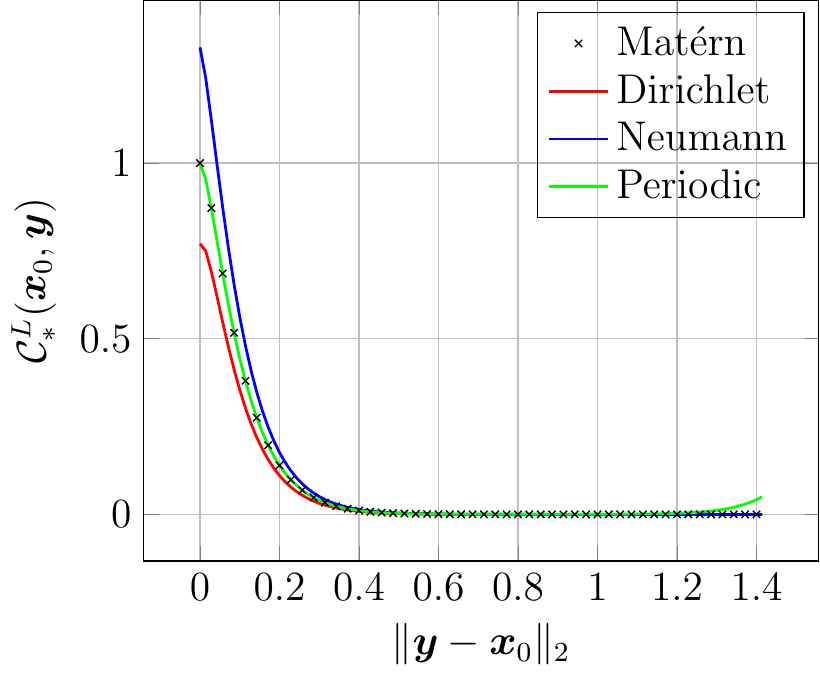}
	\caption{Case $d=2$, $\nu=1$ and~$\rho=0.1$. Covariance functions $\cov^L_{\ast}(\bfx_0,\bfy)$, $\ast\in\left\{P,D,N\right\}$, in comparison with the Mat\'ern kernel. Here $\bfx_0=(\frac{\delta}{2},\frac{\delta}{2})$ and $\bfy$ are points on the diagonal of $D$ from $(\frac{\delta}{2},\frac{\delta}{2})$ to $(1+\frac{\delta}{2},1+\frac{\delta}{2})$. Left: $\delta=0$. Center: $\delta=\rho$. Right: $\delta=2\rho$.}\label{fig:2dcovdifferentbc}
\end{figure}

\begin{figure}
	\centering
	\includegraphics[height=0.25\textwidth]{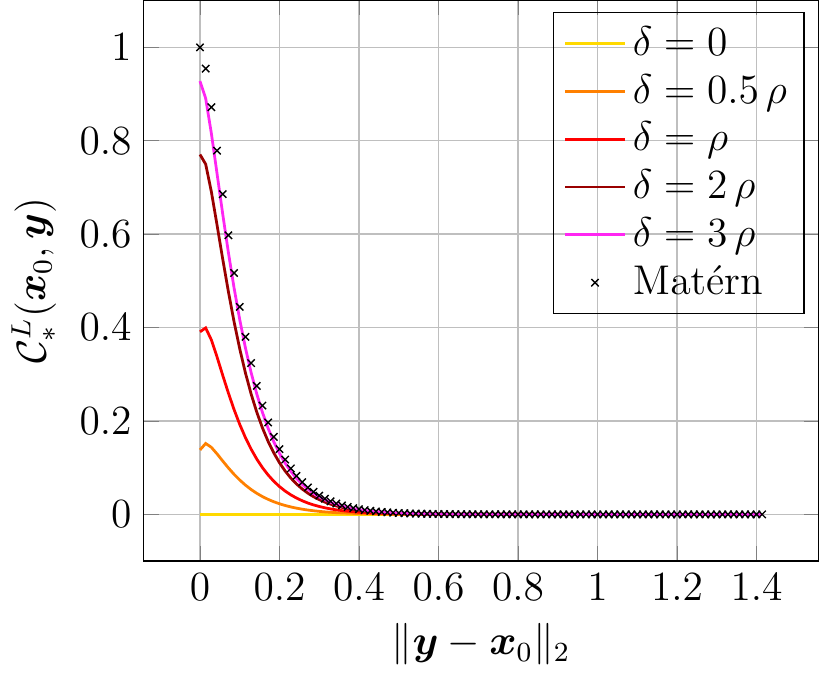}
	\hfill
	\includegraphics[height=0.25\textwidth]{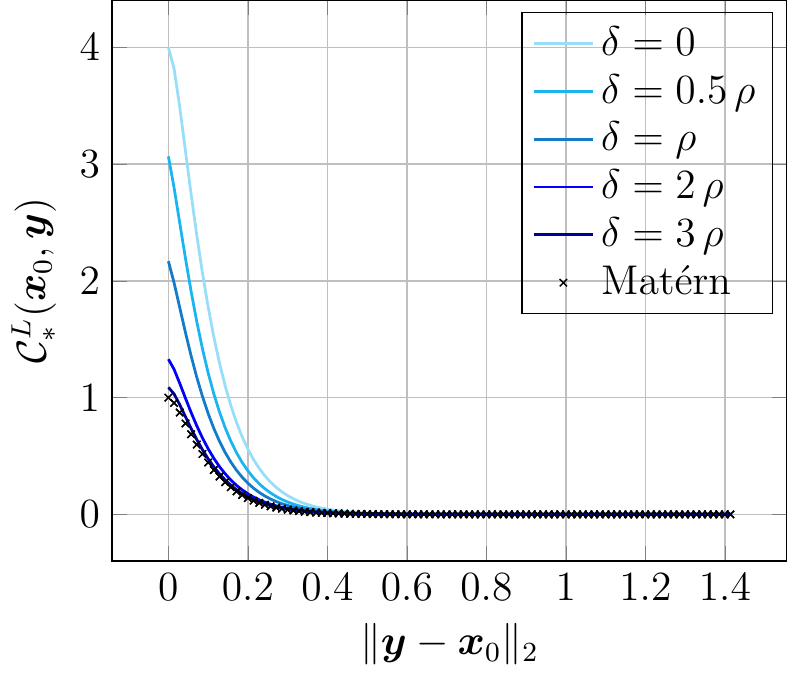}
	\hfill
	\includegraphics[height=0.25\textwidth]{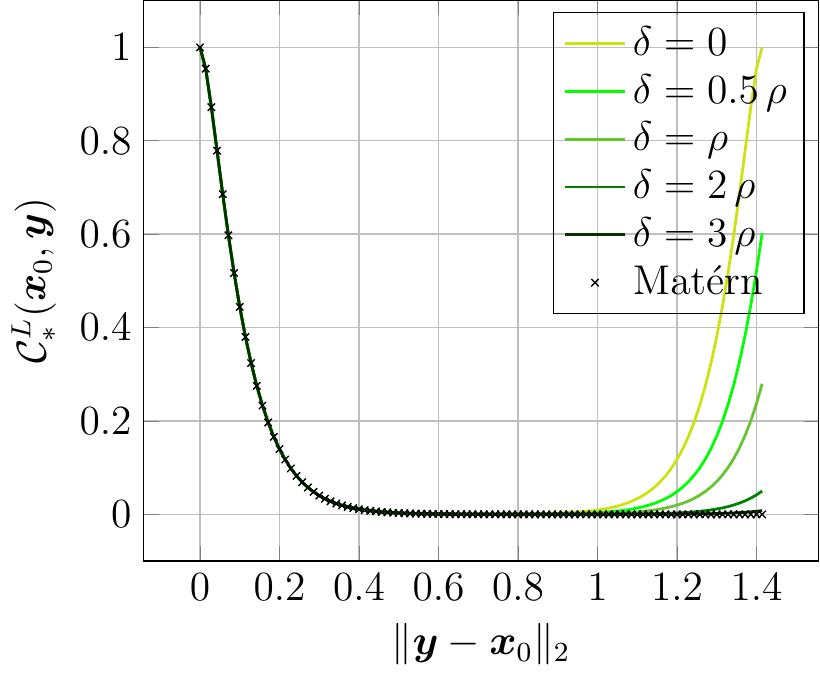}
	\caption{Case $d=2$, $\nu=1$ and $\rho=0.1$. Covariance functions $\cov^L_{\ast}(\bfx_0,\bfy)$, $\ast\in\left\{P,D,N\right\}$, in comparison with the Mat\'ern kernel, for different window sizes. Here $\bfx_0=(\frac{\delta}{2},\frac{\delta}{2})$ and $\bfy$ are points on the diagonal of $D$ from $(\frac{\delta}{2},\frac{\delta}{2})$ to $(1+\frac{\delta}{2},1+\frac{\delta}{2})$. Left: Dirichlet b.c. ($\ast=D$). Center: Neumann b.c. ($\ast=N$). Right: periodic b.c. ($\ast=P$).}\label{fig:varydelta2d}
\end{figure}

\begin{figure}
	\centering
	\begin{tabular}{ccc}
		\includegraphics[height=0.24\textwidth]{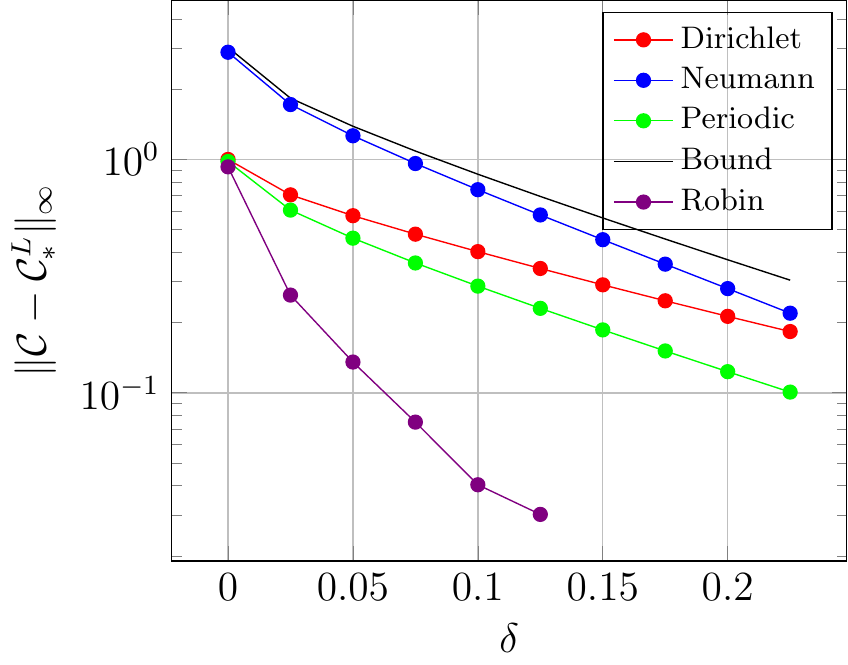} &	\includegraphics[height=0.24\textwidth]{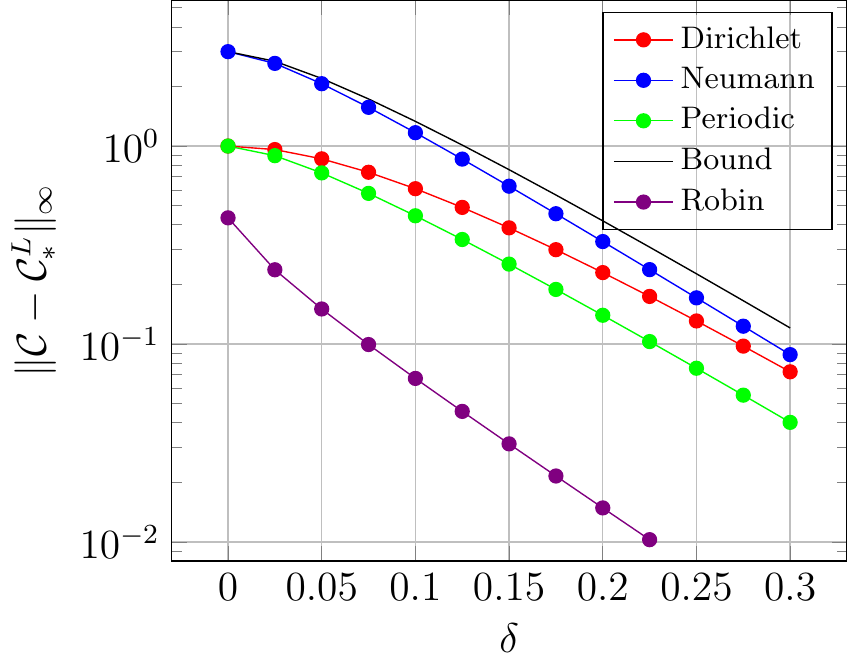} & \includegraphics[height=0.24\textwidth]{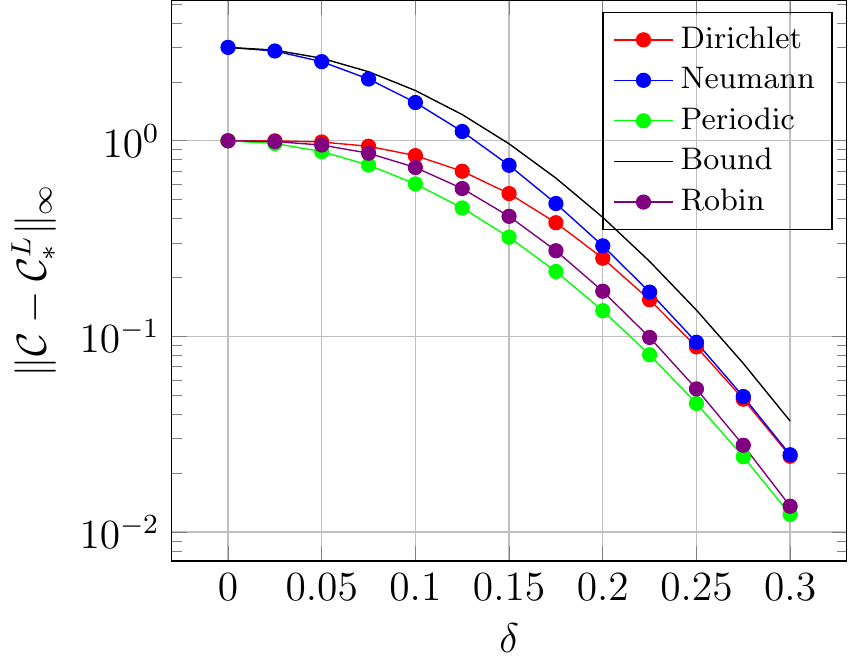}\\
		\includegraphics[height=0.24\textwidth]{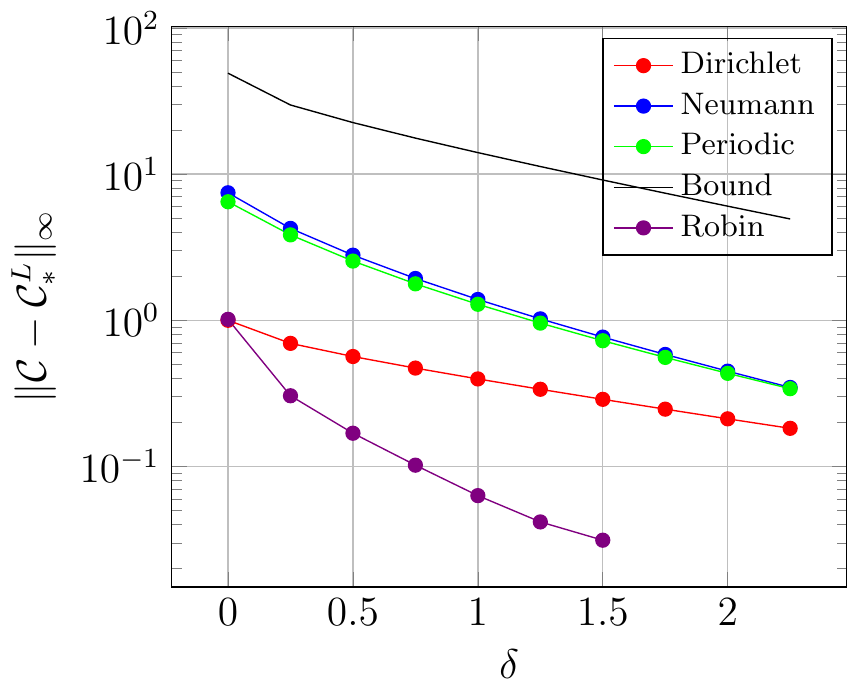} &	\includegraphics[height=0.24\textwidth]{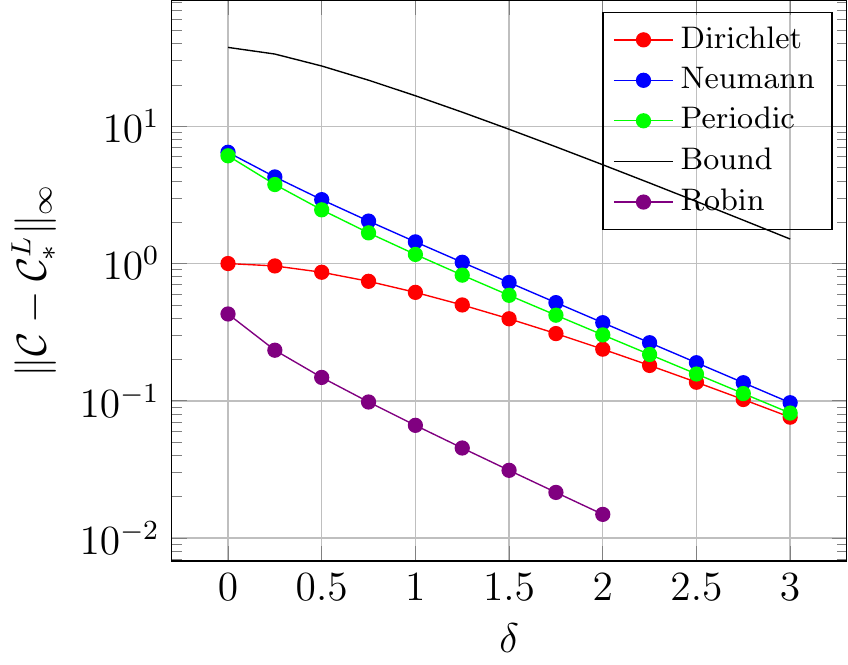} & \includegraphics[height=0.24\textwidth]{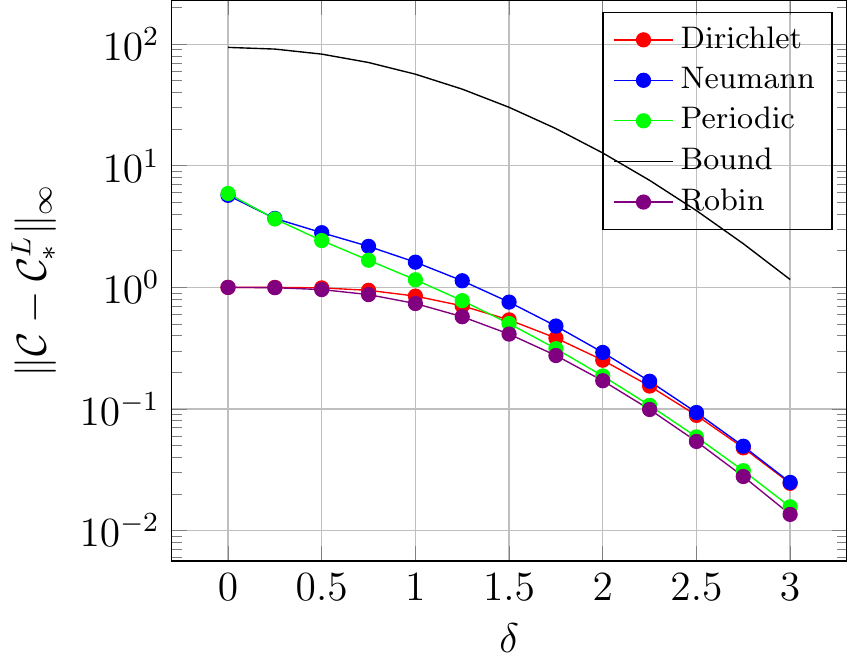}
	\end{tabular}
	\caption{Case $d=2$. Maximum norm error in the covariance kernel as function of the window size~$\delta$ for different boundary conditions. Left column: $\nu=0.25$. Center column: $\nu=1$. Right column: $\nu=50$. Top row: $\rho=0.1$. Bottom row: $\rho=1$.}\label{fig:bound2d}
\end{figure}

In these experiments, we consider $D=\left(\frac{\delta}{2},1+\frac{\delta}{2}\right)^2$, $\delta\geq 0$, as every bounded domain can be enclosed in a square or a rectangle. 
\cref{fig:2dcovdifferentbc} shows the covariance functions along the diagonal of $D$ going from $(\frac{\delta}{2},\frac{\delta}{2})$ to $(1+\frac{\delta}{2},1+\frac{\delta}{2})$, for different boundary conditions. 
We observe a similar qualitative behavior as for $d=1$ (cf. \cref{fig:1dcovdifferentbc}), but with larger errors in the approximate covariances, in particular when using Neumann boundary conditions. 
This larger error compared to the one-dimensional case can also be observed in the central plot in \cref{fig:varydelta2d} (to be compared to the central plot in \cref{fig:varydelta}). 
The latter depicts the covariances along the diagonal of $D$ for different values of the window size, each plot referring to a boundary condition. 
From \cref{fig:2dcovdifferentbc} and \cref{fig:varydelta2d}, we see that also when using Dirichlet boundary conditions the error increases with the dimension $d$, but in a much milder way than in the Neumann case. 
The periodic boundary conditions are the ones that suffer less from this effect in this case. 
Such behavior is confirmed by \cref{fig:bound2d}, reporting the error \eqref{eq:linferr} as a function of the window size. 
The error has been approximated by $\max_{\bfx,\bfy\in \mathcal{G}\times\mathcal{G}}\left|\Matern(\bfx,\bfy)-\cov^L_{\ast}(\bfx,\bfy)\right|$, for $\mathcal{G}$ a grid of $n$ equispaced points in $\left[\frac{\delta}{2},1+\frac{\delta}{2}\right]$ (extrema included). 
We have considered $n=5$ (that is $25$ points for the tensor product grid $\mathcal{G}\times\mathcal{G}$) for $\nu=1$ and $\nu=50$, and $n=3$ for $\nu=0.25$ as this case is computationally more intense (for the same reason as in the one-dimensional case).

In \cref{fig:bound2d}, we can see that the bound of \cref{thm:mainres} is sharp for Neumann boundary conditions and correlation lengths significantly smaller than the size of the domain (see the top row). 
For larger correlation lengths, the bound is not as sharp, as was observed for the one-dimensional case. 
In all the cases considered in \cref{fig:bound2d}, the Neumann boundary conditions give the largest error. 
For Dirichlet boundary conditions, the bound always overestimates the error for $\delta$ small, but for larger window sizes the error approaches eventually the same error as for Neumann boundary conditions. 
Periodic boundary conditions provide a smaller error for moderate correlation lengths, in accordance with the right plot of \cref{fig:varydelta2d}.
However, as $\rho$ increases ($\rho=1$ in the bottom row of \cref{fig:bound2d}), a similar behavior as for Neumann boundary conditions is observed. 
We can also observe, in \cref{fig:bound2d}, that Robin boundary conditions provide smaller errors than the other boundary conditions for moderate values of $\nu$, namely $\nu=0.25$ and $\nu=1$. 
When $\nu=50$, then $\kappa$ is also large and Robin boundary conditions become closer to Dirichlet boundary conditions, resulting in more similar error decay curves, see left plot in \cref{fig:bound2d}.

\begin{figure}
	\centering
	\includegraphics[height=0.27\textwidth]{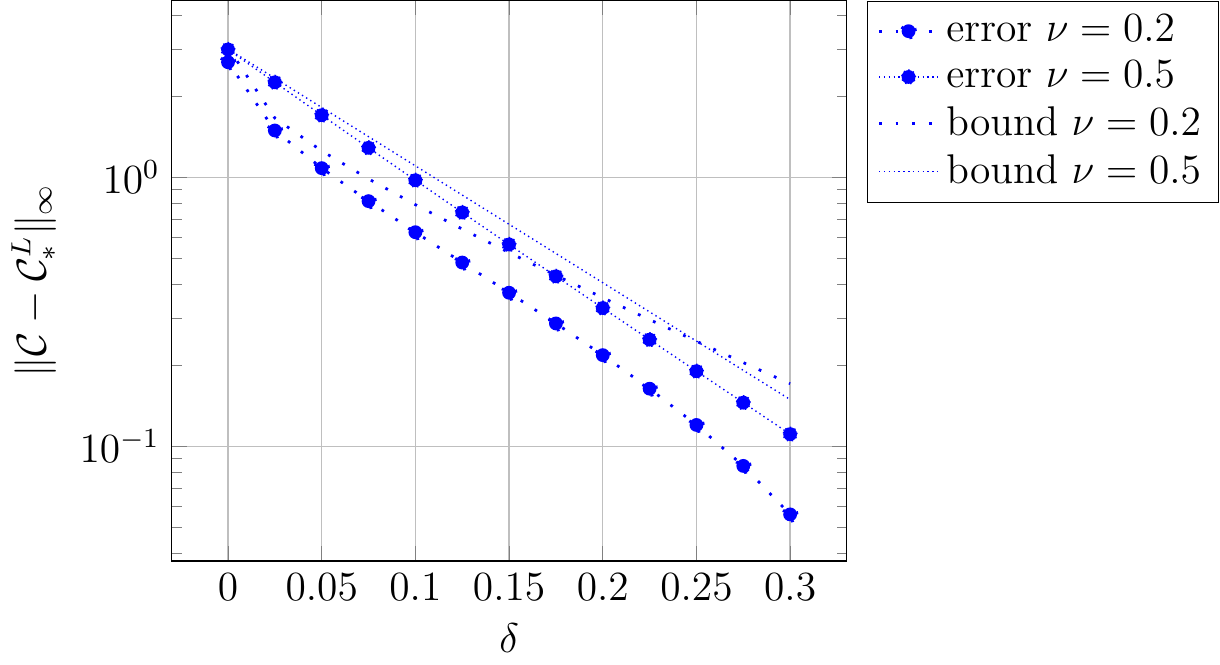}
	\hfill
	\includegraphics[height=0.27\textwidth]{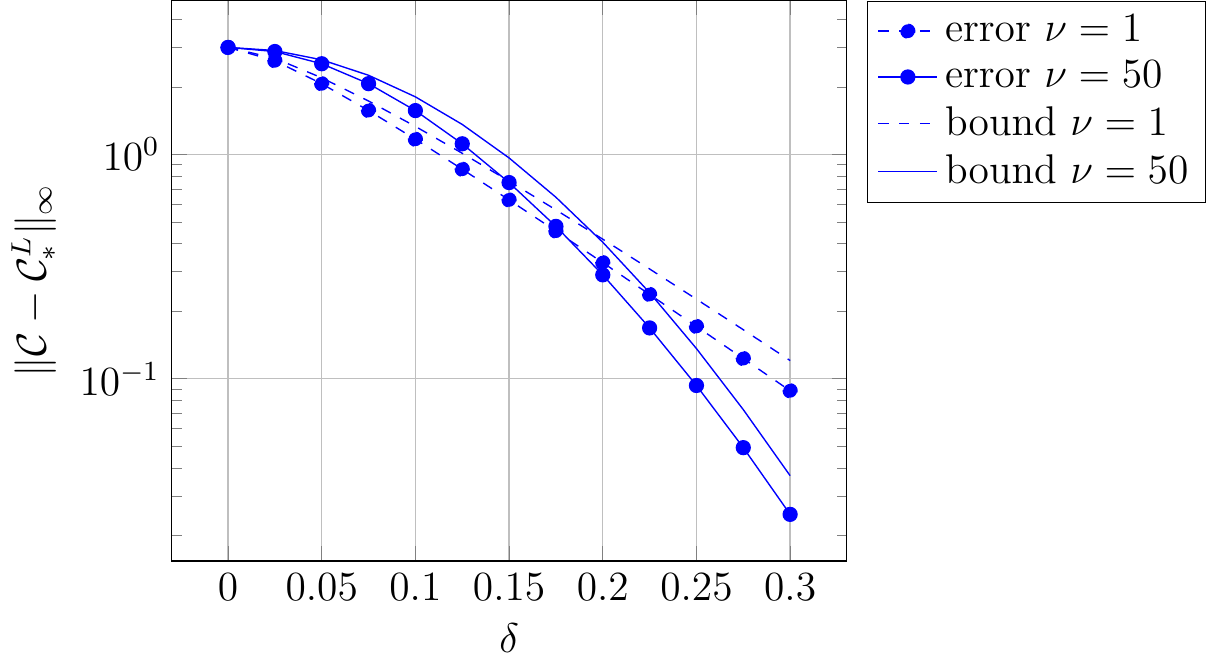}
	\caption{Case $d=2$, $\rho=0.1$, Neumann boundary conditions. Maximum norm error in the covariance kernel as function of the window size~$\delta$. Left: error and bound from \cref{thm:mainres} for different values of~$\nu\leq 0.5$.
		Right: error and bound from \cref{thm:mainres} for different values of~$\nu> 0.5$.}\label{fig:varynu2d}
\end{figure}

\begin{figure}
	\centering
	\includegraphics[height=0.3\textwidth]{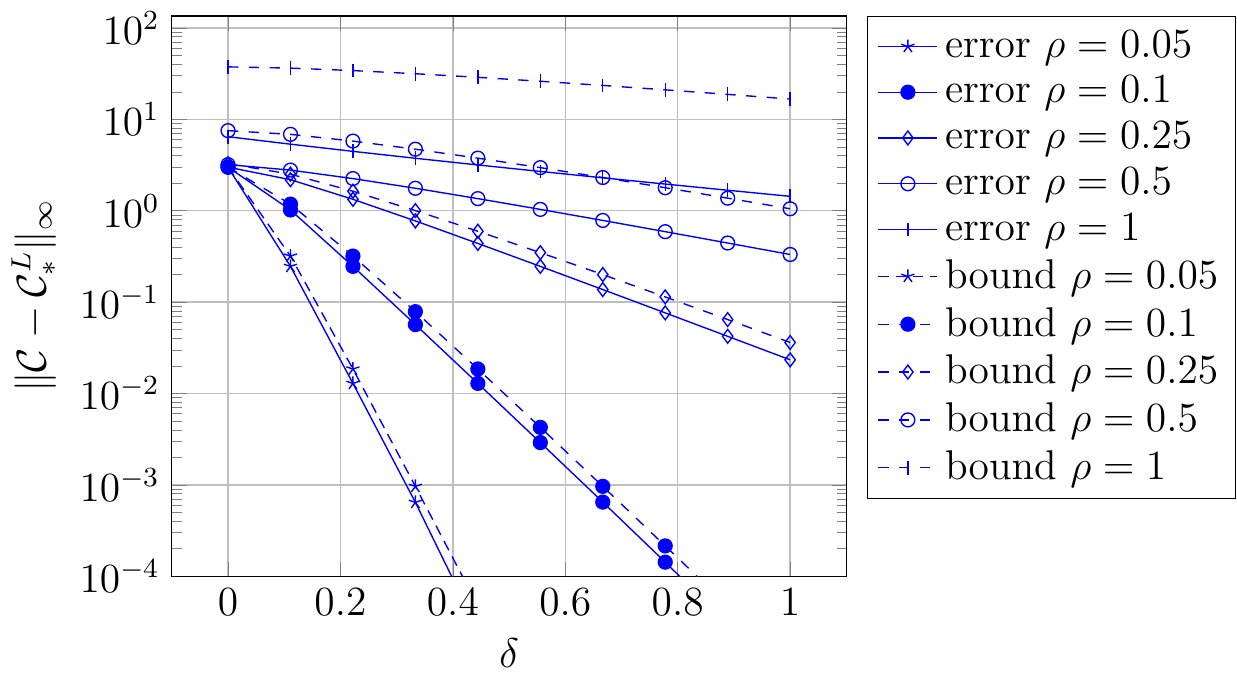}
	\caption{Case $d=2$, $\nu=1$, Neumann boundary conditions. Maximum norm error in the covariance kernel as function of the window size~$\delta$. Error and bound for different values of $\rho$.}\label{fig:varyrho2d}
\end{figure}

Finally, we consider Neumann boundary conditions and compare the behavior of the maximum norm error for different values of $\nu$ (\cref{fig:varynu2d}) and different values of $\rho$ (\cref{fig:varyrho2d}). 
Similar observations as for the one-dimensional case hold. 
In \cref{fig:varynu2d}, we see that the error bounds are less sharp than for $d=1$, but still the error behaves as $\mtn{\nu}(\sqrt{2\nu}\,\delta/\rho)$. In \cref{fig:varyrho2d}, we can see that the bound in \cref{thm:mainres} is not sharp for very large correlation lengths.
However, as already observed for the one-dimensional case, $\rho=1$ corresponds to an extreme case which is not typical in realistic simulations where the size of the domain is often larger than the correlation length. 
The fact that the error bounds are not sharp for very large $\rho$, as observed in bottom row of \cref{fig:bound1d} and \cref{fig:bound2d}, and in \cref{fig:varyrho2d}, can be explained by considering the expressions \eqref{eq:A}--\eqref{eq:ffunc} for the error bound. 
For fixed $\nu$, the value of $\kappa$ as from \eqref{eq:matcov} decreases as $\rho$ increases, leading to larger values for $f(\ell)$. 
Since $f(\ell)<1$ always, the denominator $(1-f(\ell))^d$ in \eqref{eq:A} approaches zero, and overall the value of $2^d d!f(\ell)(1-f(\ell))^{-d}$ is larger when $\rho$ is larger.



\section{Conclusions}\label{sect:concl} 
In this work we performed an error analysis for the covariance when using a window technique in PDE-based sampling of Whittle-Mat\'ern random fields. 
We considered homogeneous Dirichlet, homogeneous Neumann and periodic boundary conditions on the truncated domain. 
We have shown that, in all three cases, the error in the maximum norm decays as the Mat\'ern kernel with respect to the window size. 
Therefore, asymptotically, the error decays exponentially, with a decay rate depending on the smoothness parameter $\nu$ of the kernel. 
Numerical experiments in one and two space dimensions confirm that the error bound is sharp if the correlation length is significantly smaller than the size of the domain (in our experiments, for correlation lengths of $10\%$ the size of the domain). 
Dirichlet boundary conditions produce a smaller error than Neumann and periodic boundary conditions in the preasymptotic regime. 
In the asymptotic regime, although the error bound for the Dirichlet case is slightly sharper (cf. \eqref{eq:better}), for large $\delta$ the error curves for Dirichlet boundary conditions approach the error curves for Neumann boundary conditions, the latter providing the largest errors; periodic boundary conditions, instead, deliver a lower error than Dirichlet and Neumann boundary conditions for not too large correlations lengths.
It is evident from the theoretical results and the numerical experiments that, for $\nu$ fixed and all boundary conditions considered, the window size needed to guarantee a prescribed error in the covariance depends on the correlation length: the decay of the error depends on the argument $\kappa\delta$ of the Mat\'ern function (see \eqref{eq:bound1}), and therefore ultimately on the ratio ${\delta}/{\rho}$.
In this paper, we have also performed a numerical comparison with Robin boundary conditions with coefficient $\kappa$, and observed that for not too large values of $\nu$ these conditions deliver smaller errors in the covariance than Dirichlet, Neumann and periodic boundary conditions. 
Robin boundary conditions have the drawback that they cannot be incorporated trivially with the solution of the SPDE when using the FFT.
\eu{However,} due to their good approximation properties, they deserve further investigation. 
We foresee that the theoretical analysis of boundary effects with \eu{Robin} boundary conditions will need different techniques than those used in this paper, and we postpone it to future work.

\appendix
\section{Auxiliary results}\label{app}
\input{appendix}

\bibliographystyle{siamplain}
\bibliography{refs}

\end{document}

%% file: ex_shared.tex
\usepackage{lipsum}
\usepackage{amsfonts}
\usepackage{graphicx}
\usepackage{epstopdf}
\usepackage{algorithmic}
\ifpdf
  \DeclareGraphicsExtensions{.eps,.pdf,.png,.jpg}
\else
  \DeclareGraphicsExtensions{.eps}
\fi

\usepackage{enumitem}
\setlist[enumerate]{leftmargin=.5in}
\setlist[itemize]{leftmargin=.5in}

\newsiamremark{remark}{Remark}
\newsiamremark{hypothesis}{Hypothesis}
\crefname{hypothesis}{Hypothesis}{Hypotheses}
\newsiamthm{claim}{Claim}

\headers{Analysis of boundary effects on PDE-based sampling of Mat\'ern random fields}{U. Khristenko, L. Scarabosio, P. Swierczynski, E. Ullmann, B. Wohlmuth}

\title{Analysis of boundary effects on PDE-based sampling of Whittle-Mat\'ern random fields\thanks{Submitted to the editors DATE.
\funding{This project has received funding from the German Science Foundation, DFG grant WO-671 11-1, and from the European Union's Horizon 2020 research and innovation programme under grant agreement No 800898.}}}

\author{U. Khristenko\thanks{Technical University of Munich, Germany, Department of Mathematics, Chair of Numerical Mathematics (M2)  (\email{khristen@ma.tum.de},\email{scarabos@ma.tum.de},\email{swierczy@ma.tum.de},\email{eullmann@ma.tum.de},\email{wohlmuth@ma.tum.de}).}
\and L. Scarabosio\footnotemark[2]
\and P. Swierczynski\footnotemark[2]
\and E. Ullmann\footnotemark[2]
\and B. Wohlmuth\footnotemark[2]}

\usepackage{amsopn}


%% file: macros.tex
\newcommand{\bfeps}{\boldsymbol{\varepsilon}}

\newcommand{\bfk}{\boldsymbol{k}}
\newcommand{\bfn}{\boldsymbol{n}}
\newcommand{\bfx}{\boldsymbol{x}}
\newcommand{\bfxi}{\boldsymbol{\xi}}
\newcommand{\bfy}{\boldsymbol{y}}
\newcommand{\bfz}{\boldsymbol{z}}

\newcommand{\bfL}{\boldsymbol{L}}

\newcommand{\bbC}{ \mathbb{C} }
\newcommand{\bbE}{ \mathbb{E} }

\newcommand{\bbN}{ \mathbb{N} }

\newcommand{\bbR}{ \mathbb{R} }

\newcommand{\bbZ}{ \mathbb{Z} }

\newcommand{\dd}{ \,\mathrm{d} }

\newcommand{\img}{\textrm{i}}
\newcommand{\cov}{ \mathcal{C} }
\newcommand{\mtn}[1]{ \mathcal{M}_{#1} }
\newcommand{\Matern}{ \mathcal{C} }

\newcommand{\abs}[1]{\lVert {#1} \rVert_2}
\newcommand{\ev}{\boldsymbol{1}}
\newcommand{\Li}[1]{\operatorname{Li}_{#1}}
\newcommand{\EuNum}[2]{\left\langle\begin{matrix}
		#1\\ #2
	\end{matrix}\right\rangle}
\newcommand{\BesselK}[1]{\mathcal{K}_{#1}}

%% file: appendix.tex
In this appendix, we prove the logarithmic subadditivity of the Mat\'ern covariance function with unit marginal variance; this is used in the proof of \cref{thm:mainres}.


\begin{lemma}[Logarithmic subadditivity of unit Mat\'ern covariance]\label{lem:covbound} 
For every $x,y\in\bbR$ such that $0\leq x\leq y$, the unit Mat\'ern function $\mtn{\nu}(x)=\frac{x^{\nu}\mathcal{K}_{\nu}(x)}{2^{\nu-1}\Gamma(\nu)}$ satisfies
	\begin{equation}\label{eq:C_to_Ck}
	\begin{aligned}
	\mtn{\nu}(x+y) &\le \mtn{\nu}(x)\cdot\mtn{\nu}(y), &&\quad \text{if }\nu \ge 1/2, \\
	\mtn{\nu}(x+y) &\le \mtn{\nu}(x)\cdot\mtn{\frac{1}{2}}(y), &&\quad \text{if }0<\nu \le 1/2. \\
	\end{aligned}
	\end{equation}
\end{lemma}

\begin{proof}
We start by proving the first inequality. For this, we report here two identities \cite[\S 6.22 (15) and \S 13.71 (1)]{watson1995treatise} and one inequality \cite[(1.4)]{ismail1990complete} that are used in this proof.
For $\nu\ge 1/2$, it holds
	\begin{align}
	\BesselK{\nu}(x) &= \frac{1}{2} (x/2)^\nu\int_{0}^{\infty} e^{-t - \frac{x^2}{4t}}\frac{\dd t}{t^{\nu+1}}, \label{eq:K(2x)}\\
	\BesselK{\nu}(x)\BesselK{\nu}(y) &= \frac{1}{2} \int_{0}^{\infty} e^{-\frac{t}{2} - \frac{x^2 + y^2}{2t}}\BesselK{\nu}\left(\frac{xy}{t}\right)\frac{\dd t}{t}, \label{eq:KxK}\\	
	\frac{x^\nu\BesselK{\nu}(x)}{2^{\nu-1}\Gamma(\nu)} &\ge e^{-x}. \label{eq:K>exp}
	\end{align}	
	Using the symmetry~$\BesselK{\nu}(x)=\BesselK{-\nu}(x)$ \cite[Eq. (8) in Sect. 3.71]{watson1995treatise} and applying \eqref{eq:K(2x)} to $\BesselK{-\nu}(x)$, we arrive at
	\begin{align*}
		\mtn{\nu}(x)
		&= \frac{x^\nu}{2^{\nu-1}\Gamma(\nu)} \BesselK{-\nu}(x) = \frac{x^\nu}{2^{\nu-1}\Gamma(\nu)}\cdot\frac{1}{2^{1-\nu}} x^{-\nu}\int_{0}^{\infty} e^{-t - \frac{x^2}{4t}}\frac{\dd t}{t^{1-\nu}} \\
		&= \frac{1}{\Gamma(\nu)}\int_{0}^{\infty} t^{\nu-1}e^{-t - \frac{x^2}{4t}} \dd t.
	\end{align*}	
	Then, applying \eqref{eq:K>exp}, the change of the integration variable $t\to 2t$ and~\eqref{eq:KxK}, we obtain, for~$\nu\ge 1/2$,
	\begin{align*}
		\mtn{\nu}(x+y)
		&= \frac{1}{\Gamma(\nu)}\int_{0}^{\infty} t^{\nu-1}e^{-t - \frac{(x+y)^2}{4t}} \dd t 
		= \frac{1}{\Gamma(\nu)}\int_{0}^{\infty} t^{\nu-1}e^{-t - \frac{x^2+y^2}{4t}}e^{- \frac{xy}{2t}} \dd t \\
		&\le \frac{1}{\Gamma(\nu)}\int_{0}^{\infty} t^{\nu-1}e^{-t - \frac{x^2+y^2}{4t}}\cdot\frac{\left(\frac{xy}{2t}\right)^\nu\BesselK{\nu}\left(\frac{xy}{2t}\right)}{2^{\nu-1}\Gamma(\nu)} \dd t \\
		&\le \frac{(xy)^\nu}{(2^{\nu-1}\Gamma(\nu))^2}\cdot\frac{1}{2}\int_{0}^{\infty} e^{-t - \frac{x^2+y^2}{4t}}\cdot\BesselK{\nu}\left(\frac{xy}{2t}\right) \frac{\dd t}{t} \\
		&\le \frac{(xy)^\nu}{(2^{\nu-1}\Gamma(\nu))^2}\cdot\frac{1}{2}\int_{0}^{\infty} e^{-\frac{t}{2} - \frac{x^2+y^2}{2t}}\cdot\BesselK{\nu}\left(\frac{xy}{t}\right) \frac{\dd t}{t} \\
		&\le \frac{(xy)^\nu}{(2^{\nu-1}\Gamma(\nu))^2}\cdot\BesselK{\nu}(x)\BesselK{\nu}(y) \le \mtn{\nu}(x)\cdot\mtn{\nu}(y). \\
	\end{align*}

The second inequality in \eqref{eq:C_to_Ck} follows from the inequality~(3.2) in \cite{laforgia1991bounds},
	\begin{equation*}
	\frac{\BesselK{\nu}(x)}{\BesselK{\nu}(y)} \ge \left(\frac{y}{x}\right)^\nu e^{y-x}, \qquad 0<\nu\le \frac{1}{2}, \quad 0<x\le y,
	\end{equation*}
	that is,
\[
		\frac{\BesselK{\nu}(x)}{\BesselK{\nu}(x+y)} \ge \left(\frac{x+y}{x}\right)^{\nu}e^{y}, \qquad 0<\nu\le\frac{1}{2}, \quad x,y>0.\]
\end{proof}